\documentclass[11pt, a4paper, reqno, oneside]{article}
\usepackage{amssymb,amsmath}
\usepackage{tikz-cd}
\usepackage{authblk} % For affiliations
\usepackage{physics} % For derivative notation (and more)
\usepackage{tikz}
\usepackage{titlesec} % Section/subsection/subsubsection style
\usetikzlibrary{babel}
\usepackage{adjustbox}
\usepackage{url}
\usepackage{xcolor}
\usepackage{eso-pic}
\usepackage[T1]{fontenc}
\usepackage[utf8]{inputenc}
\usepackage{csquotes}
\usepackage{chngcntr}
\usepackage{apptools}
\usepackage{pst-node}
\usepackage[english]{babel}
\usepackage{libertine}
\usepackage[libertine]{newtxmath}
\usepackage{latexsym}
\usepackage{euscript}
\usepackage{mathtools}
\usepackage[all]{xy}
\usepackage{mathrsfs}
\usepackage[percent]{overpic}
\usepackage{setspace}
\usepackage{graphicx}
\usepackage{pdfsync}

\usepackage{amsthm}
\usepackage{thmtools}
\usepackage{yfonts}

\usepackage[colorlinks=true,linktocpage=true, pagebackref=false, citecolor=black,linkcolor=black]{hyperref}
\usepackage{cite}

\normalbaselines
\makeatletter

\def\ps@myfancy{\let\@mkboth\markboth
 \def\@evenhead{\vbox{\hsize\textwidth 
 \hbox to \textwidth{\sf\mdseries\thepage 
 \rule[-.6ex]{0mm}{2mm} \hfill\sf\large\leftmark}
 \vskip 1pt \hrule}}
 \def\@oddhead{\vbox{\hsize\textwidth 
 \hbox to \textwidth{{\sf\large\leftmark}
 \rule[-.6ex]{0mm}{2mm} \hfill\sf\mdseries{\thepage}}
 \vskip 1pt \hrule}}}

\def\ps@myfancyplain{
 \def\@evenhead{\vbox{\hsize\textwidth%
 \rule[-.6ex]{0mm}{2mm} \hfill }
 \vskip 1pt \hrule
 \vskip\headsep
 \vskip\textheight
 \vskip1pc
 \hbox to \textwidth{\sf\mdseries\thepage 
 \rule[-6ex]{0mm}{2mm} \hfill }}
 \def\@oddhead{\vbox{\hsize\textwidth 
 \vskip 1pt\hrule
 \vskip\headsep
 \vskip\textheight
 \vskip2pc
 \hbox to \textwidth{\hfill\rule[.4ex]{1pc}{2.5pt}
 \sf\mdseries\thepage}
}}}

\def\ps@myemptyfun{
 \def\@evenhead{\vbox{\hsize\textwidth
 \rule[-.6ex]{0mm}{2mm} \hfill }
 \vskip 1pt 
 \vskip\headsep
 \vskip\textheight
 \vskip1pc
 \hbox to \textwidth{\sf\mdseries\thepage 
 \rule[-0.6ex]{0mm}{2mm} 
 \hfill }}
 \def\@oddhead{\vbox{\hsize\textwidth 
 \vskip 1pt 
 \vskip\headsep
 \vskip\textheight
 \vskip2pc
 \hbox to \textwidth{\hfill\rule[.5ex]{3pc}{0.8pt}
 \ \sf\mdseries\textcolor{blue}{\thepage}}
}}}

\makeatother
\providecommand{\proofname}{Demostraci\'on.}
 {\par\noindent{\it Demostraci\'on. }\nopagebreak\normalsize}%

 {\par\noindent{\it #1. }\nopagebreak\normalsize}%
 {\hfill\linebreak[2]\hspace*{\fill}$\square$\\[-1pt]}
\makeatother

\def\sqbullet{\raise.2ex\hbox{\vrule width 3.5pt height 3.5pt}}

{}

{\em}

\newcounter{substep}
\def\thesubstep{\arabic{substep}}

{\em}

\newcounter{subsubstep}
\def\thesubsubstep{\arabic{subsubstep}}

{\em}

\numberwithin{figure}{section}

{}

\newtheoremstyle{mystyle}%                % Name
  {}%                                     % Space above
  {}%                                     % Space below
  {\itshape}%                                     % Body font
  {}%                                     % Indent amount
  {\sf \bfseries}%                            % Theorem head font
  {}%                                    % Punctuation after theorem head
{ }%                                    % Space after theorem head, ' ', or \newline
  {\thmname{#1}\thmnumber{{\textcolor{blue}{\, \hspace{-1mm}#2.}}}\thmnote{ (#3)}}%                                     % Theorem head spec (can be left empty, meaning `normal')

\theoremstyle{mystyle}
\definecolor{royalblue(web)}{rgb}{0.25, 0.41, 0.88}
\hypersetup{
    colorlinks=true,
    linkcolor=royalblue(web),
    filecolor=magenta,      
    urlcolor=cyan,
    pdftitle={Overleaf Example},
    pdfpagemode=FullScreen,
    }

\titleformat{\section}
  {\normalfont\LARGE\bfseries}{\thesection.}{1em}{}
\titleformat{\subsection}
  {\normalfont\Large\bfseries}{\thesubsection.}{1em}{}
\titleformat{\subsubsection}
  {\normalfont\Large\bfseries}{\thesubsubsection.}{1em}{}

\newtheorem{Teor}{Theorem}[section]

\newtheorem{Prop}[Teor]{Proposition}
\newtheorem{Coro}[Teor]{Corollary}
\newtheorem{Lema}[Teor]{Lemma}
\newtheorem{Defi}[Teor]{Definition}

 \newcommand{\R}{{\mathbb R}}
 \newcommand{\C}{{\mathbb C}}

\newcommand{\mail}[1]{\small\href{mailto:#1}{#1}}

\usepackage{stackengine}
\newcommand\dhookrightarrow{\mathrel{%
  \ensurestackMath{\stackanchor[.1ex]{\hookrightarrow}{\hookrightarrow}}
}}

\newenvironment{keywords}
{
\begin{center}
\textbf{Keywords}\\
\vspace{0.17cm}
\begin{minipage}{14.5cm}}
{\footnotesize
\end{minipage}
\end{center}}

\newenvironment{Abstract}
{
\begin{center}
\textbf{Abstract}\\
\vspace{0.25cm}
\begin{minipage}{14.5cm}}
{\footnotesize
\end{minipage}
\end{center}}

\begin{document}
%-------------------------------------------
% Title, authors and affiliations
%-------------------------------------------
%\maketitle

	\begin{center}
		{\huge {\bfseries Compact embeddings of Bessel Potential Spaces}\par}
		\vspace{1cm}
		
\begin{tabular}{l@{\hskip 2cm}l} 
	{\Large Jos\'e C. Bellido }{\small\textsuperscript{1}} & {\Large Javier Cueto }{\small\textsuperscript{2}} \\
	\mail{josecarlos.bellido@uclm.es} & \mail{javier.cueto@uam.es} \\
	[10pt] % Salto de l\'inea gordito
\end{tabular}\\
{\centering{\Large Guillermo Garc\'ia-S\'aez}{\small\textsuperscript{1}} \\
	\mail{guillermo.garciasaez@uclm.es}} \\ \vspace{5mm}

\textsc{\textsuperscript{1}ETSI Industrial\\ Departamento de Matem\'aticas, Universidad de Castilla-La Mancha} \\
		Campus Universitario s/n, 13071-Ciudad Real, SPAIN. \\ \vspace{5mm}

\textsc{\textsuperscript{2}Departamento de Matem\'aticas, Universidad Aut\'onoma de Madrid} \\
      Campus de Cantoblanco, 28049 Madrid, SPAIN. \\
        \vspace{5mm}
		
\end{center}

\date{\today}
%-------------------------------------------
% Abstract
%-------------------------------------------
\begin{Abstract}

Bessel potential spaces have gained renewed interest due to their robust structural properties and applications in fractional partial differential equations (PDEs). These spaces, derived through complex interpolation between Lebesgue and Sobolev spaces, are closely related to the Riesz fractional gradient. Recent studies have demonstrated continuous and compact embeddings of Bessel potential spaces into Lebesgue spaces. This paper extends these findings by addressing the compactness of continuous embeddings from the perspective of abstract interpolation theory. We present three distinct proofs, leveraging compactness results, translation estimates, and the relationship between Gagliardo and Bessel spaces. Our results provide a deeper understanding of the functional analytic properties of Bessel potential spaces and their applications in fractional PDEs.

\end{Abstract}

%-------------------------------------------
% Keywords
%-------------------------------------------

%-------------------------------------------
% Table of contents
%-------------------------------------------
\tableofcontents
%\tableofcontents

%-------------------------------------------
% Introduction to the results
%-------------------------------------------

\section{Introduction}

Bessel potential spaces, introduced in the 1960s \cite{Aronszajn1961}, emerged from the study of partial differential equations (PDEs), particularly in the linear case. Notably, these spaces are the images of trace operators acting on Sobolev spaces \cite{LionsMagenes1972}. They can be derived through complex interpolation between Lebesgue and Sobolev spaces \cite{Calderon1963,Calderon1961,Lions1960,LionsMagenes1972}, situating them as intermediate spaces or Sobolev spaces of fractional differentiability order. Another class of fractional Sobolev spaces, known as Gagliardo or Sobolev-Slobodeckij spaces, is obtained via the real method of interpolation rather than the complex method. In recent decades, Gagliardo spaces have gained significant attention due to their direct connection to the fractional Laplacian, or more generally, the fractional $p$-Laplacian. Conversely, Bessel potential spaces have only recently begun to receive more attention.

There has been a renewed interest in Bessel potential spaces following the pioneering works \cite{ShiehSpector2015,ShiehSpector2018}, which established a connection with the Riesz fractional gradient. The Riesz fractional gradient of order $s$ for a smooth, compactly supported function $u$ is defined as
\[D^s u(x) = c_{n,s} \int_{\mathbb{R}^n} \frac{u(x) - u(y)}{|x - y|^{n+s}} \frac{x - y}{|x - y|} \, dy,\]
where $c_{n,s}$ is a normalization constant dependent on the dimension $n$ and the fractional index $s$.

In \cite{ShiehSpector2015,ShiehSpector2018}, it was demonstrated that Bessel potential spaces coincide, with equivalence of norms, with the closure of $C_c^\infty(\mathbb{R}^n)$ with respect to the norm
\[ \|u\| = \|u\|_{L^p(\mathbb{R}^n)} + \|D^s u\|_{L^p(\mathbb{R}^n, \mathbb{R}^n)}. \]

Recent applications of Bessel potential spaces include the study of a new family of fractional PDEs \cite{ShiehSpector2015,ShiehSpector2018}, variational inequalities involving the Riesz fractional gradient \cite{Campos2024,campos2023Primero,CamposRodrigues2023II,LoRodriges2023}, fractional elasticity problems \cite{BellidoCuetoMoraCorral2021,BellidoCuetoMoraCorral2020}, and the development of a new concept of fractional perimeter \cite{BrueCalziComiStefani2022,ComiStefani2023,ComiStefani2019}.

The increasing interest in Bessel potential spaces can be attributed to their strong structural properties from a Functional Analysis perspective, which makes them particularly valuable in the study of fractional partial differential equations (PDEs). Notably, \cite{ShiehSpector2015} demonstrated continuous embeddings of Bessel potential spaces into Lebesgue spaces, while \cite{ShiehSpector2018} established the compactness of these embeddings. These findings are fundamentally based on the identification of these spaces through the Riesz fractional gradient, as previously discussed. However, the role of Bessel potential spaces as interpolation spaces has not been fully leveraged in these analyses.

This paper is a follow-up work of \cite{BellidoGarciaA2025}, where continuous embeddings for Bessel potential spaces are studied from the complex interpolation method perspective. In this paper we address the compactness of the continuous embeddings established in \cite{BellidoGarciaA2025}. Given two Banach spaces $E,F$, we denote by $$E\dhookrightarrow F,$$ the compact embedding of $E$ into $F$. With this notation, the main result in this paper is the following:

\begin{Teor}[Compact embeddings for Bessel potential spaces]\label{BesselCompact}
    Let $s\in (0,1)$, $1\le p<\infty$ and $\Omega\subset \R^n$ be a bounded, Lipschitz, open set. Then, 
    \begin{itemize}
        \item[a)] {\bf (Subcritical case)} If $sp<n$, then $H^{s,p}(\R^n)\dhookrightarrow L^q(\Omega)$ for $1\leq q<p_s^*=\frac{np}{n-sp}$.
        \item[b)] {\bf (Critical case)} If $sp=n$, then $H^{s,p}(\R^n)\dhookrightarrow L^q(\Omega)$ for $1\leq q<\infty$.
        \item[c)] {\bf (Supercritical case)} If $sp>n$, then $H^{s,p}(\R^n)\dhookrightarrow C^{0,\mu}(\overline{\Omega})$ for $0<\mu<\mu_s^*=s-n/p$.
    \end{itemize}
\end{Teor}

This result extends the existing findings in the literature. For instance, for $1 < p < \infty$, part $a)$, which is a Rellich-Kondrachov compactness result for Bessel potential spaces, was established in \cite[Theorem 2.2]{ShiehSpector2018} under the condition of imposing function values on the complement of the domain $\Omega$. Under the same hypothesis, parts $b)$ and $c)$ were demonstrated in \cite{BellidoCuetoMoraCorral2020}. We present three distinct proofs of this result. The first proof utilizes compactness results from abstract interpolation theory, an approach that, to the best of the authors' knowledge, has not been previously explored in the literature. The second proof relies on an estimation of translations based on the Riesz fractional gradient and the application of the Fréchet-Kolmogorov theorem. The third proof leverages the well-known relationship between Gagliardo and Bessel spaces, along with compactness results for the former.

The structure of this paper is as follows. Section 2.1 provides a concise overview of the main results on compactness in abstract interpolation theory, covering both the real and complex methods. Section 2.2 summarizes the key findings from \cite{BellidoGarciaA2025} and introduces Bessel spaces on domains. Section 2.3 properly introduces the Riesz fractional gradient and presents the fractional fundamental theorem of calculus. Section 3 is dedicated to the various proofs of Theorem \ref{BesselCompact} as mentioned earlier.

\section{Preliminaries}
\subsection{Interpolation theory and compactness}
In this section we review the long-standing problem of the compactness and the theory of interpolation. We first recall the very basic definitions of interpolation theory and the two main methods: the real and the complex methods.\\

Let $(E_0,E_1)$ be a couple of Banach spaces. We say that the couple is \textit{compatible} if there exists a Hausdorff topological vector space $\mathcal{E}$ such that $E_0,E_1\xhookrightarrow{}\mathcal{E}$. We say that a Banach space $E$ is intermediate with respect to the couple if $E_0\cap E_1\xhookrightarrow{}E\xhookrightarrow{}E_0+E_1$. Let $(F_0,F_1)$ be another compatible couple and $T:E_0+E_1\to F_0+F_1$ a bounded linear operator. We would say that such a $T$ is \textit{admissible} if $T:E_j\to F_j$, $j=0,1,$ continuously. Given $E$ an intermediate space with respect to the couple $(E_0,E_1)$, and $F$ an intermediate space with respect to $(F_0,F_1)$, we say that $E,F$ are \textit{interpolation spaces} with respect to the couples $(E_0,E_1)$ and $(F_0,F_1)$, respectively, if for every admissible operator $T:E_0+E_1\to F_0+F_1$, $T:E\to F$ continuously. The methods to construct such interpolation spaces for given couples are called \textit{interpolation functors}, i.e., given a compatible couple $(A,B)$, an interpolation functor $\mathcal{F}$ maps $(A,B)$ into an interpolation space $\mathcal{F}\left((A,B)\right)$. We say that an interpolation functor $\mathcal{F}$ is of \textit{exponent} $\theta\in (0,1)$ if $$\norm{T}_{\mathcal{F}\left((E_0,E_1)\right)\to \mathcal{F}\left((F_0,F_1)\right)}\leq C\norm{T}_{E_0\to F_0}^{1-\theta}\norm{T}_{E_1\to F_1}^\theta,$$ for some positive constant $C$. If we can choose $C=1$, we say that the functor is \textit{exact} of exponent $\theta$.

The real interpolation method is defined in the following way. Let $\theta\in  (0,1)$ and $1\leq q\leq \infty$. Given a compatible couple of Banach spaces $(E_0,E_1)$, for $t>0$ and $x\in E_0+E_1$, we define the $K$-functional as $$K\left((t,x;E_0,E_1)\right):=\operatorname{inf}_{x=x_0+x_1}\{\norm{x_0}_{E_0}+t\norm{x_1}_{E_1}:x_j\in E_j,j=0,1\}.$$
If there is no confusion about the couple, we denote it simply as $K(t,x)$. We define the \textit{real interpolation space} $(E_0,E_1)_{\theta,q;K}$ as $$(E_0,E_1)_{\theta,q;K}:=\{x\in E_0+E_1: \norm{x}_{\theta,q;K}<\infty\},$$ where $$
    \norm{x}_{\theta,q;K}:=\begin{cases}
        \left(\int_0^\infty \left(t^{-\theta}K(t,x)\right)^q\,\frac{dt}{t}\right)^{1/q},\,&q<\infty,\\
        \operatorname{sup}_{t>0}t^{-\theta}K(t,x),\,&q=\infty.
    \end{cases}$$
The spaces $(E_0,E_1)_{\theta,q;K}$ are interpolation spaces (see \cite[3.1.2]{BerghLofstrom1976} or \cite[Theorem~III.1.3]{GarciaSaez2024}). Moreover, the functor $K_{\theta,q}:(E_0,E_1)\mapsto (E_0,E_1)_{\theta,q;K}$ is of exponent $\theta$.

Since we will be dealing with Bessel potential spaces, we are mainly interested in the complex method. Given a compatible couple of Banach spaces $(E_0,E_1)$, we define the space $\mathfrak{F}(E_0,E_1)$, as the space of functions $f:S\to E_0+E_1$, where $S:=\{z\in \C: 0\leq \operatorname{Re}z\leq 1\}$, such that $f$ is holomorphic on the interior of $S$, continuous and bounded on $S$, and the functions $t\mapsto f(j+it)$, $j=0,1$, are continuous from $\R\to E_j$, and such that $\norm{f(j+it)}_{E_j}\to 0$ as $|t|\to \infty$. The space $\mathfrak{F}(E_0,E_1)$ is a vector space which is complete endowed with the norm $$\norm{f}_{\mathfrak{F}(\overline{E})}:=\operatorname{max}\{\operatorname{sup}_{t\in\R}\norm{f(it)}_{E_0},\operatorname{sup}_{t\in\R}\norm{f(1+it)}_{E_1}\},\,f\in \mathfrak{F}(\overline{E}).$$ From this space, we construct the \textit{complex method} as the functor $\mathcal{C}_\theta$, $\theta\in [0,1]$ which associates the space $[E_0,E_1]_\theta$ to the compatible couple of Banach spaces $(E_0,E_1)$. The space $[E_0,E_1]_\theta$ is defined as the space of $x\in E_0+E_1$ such that there exists $f\in \mathfrak{F}(E_0,E_1)$ with $f(\theta)=x$. The space is a Banach space endowed with the norm $$\norm{x}_{\theta}:=\operatorname{inf}\{\norm{f}_{\mathfrak{F}(\overline{E})}: f\in \mathfrak{F}(\overline{E}), f(\theta)=x\}.$$ 
The spaces $[E_0,E_1]_\theta$ are exact interpolation spaces of exponent $\theta$ (see \cite[Theorem~4.1.2]{BerghLofstrom1976} or \cite[Theorem~IV.1.5]{GarciaSaez2024}).

From the very beginning of the theory, the question of whether the interpolation methods preserve the compactness of the admissible operators was posed, motivated by the result of Krasnoselskii \cite{Krasnoselskii1960}, which generalized the Riesz-Thorin theorem under the assumption of compactness in one of the endpoints, or spaces of the compatible couple. Given $(U,\mu)$ a finite positive measure, and $1\leq p\leq \infty$, we define $L^p(U,d\mu)$ as the space of $\mu$-measurable functions $f$ on $U$ such that $$\norm{f}_p:=\left(\int_U|f|^p\,d\mu\right)^{1/p}<\infty,$$ when $1\leq p<\infty$, and the corresponding one when $p=\infty$. For such spaces we have the following result:
\begin{Teor}[Krasnoselskii's theorem]\label{Krasno}
    Let $(U,\mu),(V,\nu)$ be two finite positive measure spaces and $1\leq p_0,p_1,q_1\leq\infty$, $1\leq q_0<\infty$. Assume that $$T:\left(L^{p_0}(U,d\mu),L^{q_0}(V,d\nu)\right)\to \left(L^{p_1}(U,d\mu),L^{q_1}(V,d\nu)\right),$$ is an admissible operator such that $T:L^{p_0}\to L^{q_0}$ is compact. Then, $T:L^p(U,d\mu)\to L^q(V,d\nu)$ is compact for any $$\frac{1}{p}=\frac{1-\theta}{p_0}+\frac{\theta}{p_1},\,\frac{1}{q}=\frac{1-\theta}{q_0}+\frac{\theta}{q_1},\,\theta\in (0,1).$$
\end{Teor}
There were earlier works on the subject, like the works of Kantorovich on compactness of integral operators \cite{Kantorovic1956}. These results and the foundation of abstract interpolation methods led to the question of replacing the couples of Lebesgue spaces by abstract compatible couples of Banach spaces.  One would expect that given two compatible couples of Banach spaces $\overline{E}=(E_0,E_1),\, \overline{F}=(F_0,F_1)$, and a compatible operator $T:\overline{E}\to \overline{F}$, if $T:E_i\to F_i,$ $i=0,1$ is compact, then $T:\mathcal{F}(\overline{E})\to\mathcal{F}(\overline{F})$ is compact, given some interpolation functor $\mathcal{F}$. These are the two-sided compactness results, since we are imposing that the operator $T$ is compact on both endpoints. The stronger version is the one-sided one, which only assumes that $T$ is compact acting from $E_0\to F_0$ or $E_1\to F_1$. We will first review the problem from the real interpolation point of view and then for the complex method. \\
For the real method, the first results were obtained by Lions and Peetre \cite{LionsPeetre1964} (see also the lectures of Gagliardo \cite{Gagliardo1961}). They established two one-sided results for degenerate couples, i.e., couples $(E_0,E_1)$ such that $E_0=E_1$, and one result concerning ordered couples ($E_0\xhookrightarrow{}E_1)$ with compact inclusion, that we summarize in the following theorem. 
%\begin{Teor}[Lions-Peetre First Compactness Theorem]\label{lpI}
 %   Let $\overline{F}=(F_0,F_1)$ a compatible couple of Banach spaces, $F\in \Cc_J\left(\theta;\overline{F}\right)$, $0<\theta<1$, and $E$ a Banach space. Suppose that $T\in \mathcal{L}(E,F_i)$, $i=0,1$, and $T:E\to F_0$ compactly. Then, $T:E\to F$ compactly.
%\end{Teor}
%\begin{Teor}[Lions-Peetre Second Compactness Theorem]\label{lpII}
    %Let $\overline{E}=(E_0,E_1)$ a compatible couple of Banach spaces, $E\in \Cc_K\left(\theta;\overline{E}\right)$, $0<\theta<1$, and $F$ a Banach space. Suppose that $T\in \mathcal{L}(E_i,F)$, $i=0,1$, and $T:E_0\to F$ compactly. Then, $T:E\to F$ compactly.
%\end{Teor}
%\begin{Teor}[Lions-Peetre Third Compactness Theorem]
    %Let $(E_0,E_1)$ a compatible couple of Banach spaces, $0<\theta_0<\theta_1<1$ and $F_i\in \Cc\left(\theta_i;(E_0,E_1)\right)$, $i=0,1$. If $E_0\xhookrightarrow{}E_1$ compactly, then $F_0\xhookrightarrow{}F_1$ compactly.
%\end{Teor}
%From here the following compactness result for the real method is straightforward. 
\begin{Teor}[Lions-Peetre]
    Let $(E_0,E_1), (F_0,F_1)$ be compatible couples of Banach spaces. Let $0<\theta_0,\theta,\theta_1<1$ such that $\theta_0<\theta_1$, and $1\leq q_0,q,q_1\leq \infty$ Let $T:(E_0,E_1)\to (F_0,F_1)$ an admissible operator such that $T:E_0\to F_0$ compactly.
    \begin{itemize}
        \item If $E_0=E=E_1$, $T:E\to (F_0,F_1)_{\theta,q}$ compactly.
        \item If $F_0=F=F_1$, $T:(E_0,E_1)_{\theta,q}\to F$ compactly.
    \end{itemize}
    Furthermore, if $E_0\xhookrightarrow{}E_1$ compactly, then $(E_0,E_1)_{\theta_0,q_0}\xhookrightarrow{} (E_0,E_1)_{\theta_1,q_1}$ compactly.
\end{Teor}
Additionally, since for every interpolation functor of exponent $\theta\in (0,1)$ produces spaces of $\theta$-class, we can generalize previous results for more general functors. \begin{Teor}[Functorial first and second Lions-Petree Theorem]\label{FLP}
   Let $\overline{E}=(E_0,E_1), \overline{F}=(F_0,F_1)$ compatible couples of Banach spaces, $\mathcal{F}$ an interpolation functor of exponent $\theta\in (0,1)$, and $T:(E_0,E_1)\to (F_0,F_1)$ an admissible operator such that $T:E_0\to F_0$ compact.
    \begin{itemize}
        \item If $E_0=E=E_1$, $T:E\to \mathcal{F}(\overline{F})$ compactly.
        \item If $F_0=F=F_1$, $T:\mathcal{F}(\overline{E})\to F$ compactly.
    \end{itemize}
\end{Teor}
\begin{Coro}[Functorial Lions-Peetre Third Theorem]\label{functorialPeetreIII}
    Let $\overline{E}=(E_0,E_1)$ be a compatible couple of Banach spaces such that $E_0$ is compactly embedded in $E_1$. If $\mathcal{F}$ and $\mathcal{G}$ are interpolation functors of exponents $\theta_0$ and $\theta_1$, respectively, with $0<\theta_0<\theta_1<1$, then $\mathcal{F}(\overline{E})$ is compactly embedded into $\mathcal{G}(\overline{E})$.
\end{Coro}
%Soon, also in 1964, Persson \cite{Persson1964} proved the first compactness result for non-degenerate couples, assuming the following approximation hypothesis on the second couple: Let $(F_0,F_1)$ be a compatible couple of Banach spaces such that for every compact subset $K\subset E_0$, there exists a constant $C>0$ and a set $\{T_\alpha\}_{\alpha\in \Lambda}$ of linear operators $T_\alpha:F_0+F_1\to F_0\cap F_1$ such that $$\norm{T_\alpha}_{\mathcal{L}(E_j,E_j)}\leq C,\,j=0,1,\,\alpha\in \Lambda,$$ and such that for every fixed $\varepsilon>0$ there exists $T_\alpha$ such that $$\norm{T_\alpha x-x}_{E_0}<\varepsilon,\,x\in K.$$
%The approximation hypothesis corresponds to Kranoselskii's assumption of $q_0<\infty$ in \textbf{Theorem \ref{Krasno}} (see \cite{Persson1964} and \cite[Lemma~2.4]{Cobos2009}).
%\begin{Teor}[Persson]\label{Persson}
 %   Let $(E_0,E_1),(F_0,F_1)$ two compatible couples of Banach spaces, and $E,F$ two interpolation spaces of exponent $\theta\in (0,1)$ with respect to the couples $(E_0,E_1), (F_0,F_1)$, respectively. Suppose that $E\in \mathcal{C}_{K}\left(\theta;(E_0,E_1)\right)$ and that $(F_0,F_1)$ satisfy the approximation hypothesis. Then, if $T:(E_0,E_1)\to (F_0,F_1)$ is an admissible operator such that $T:E_0\to F_0$ is compact, then $T:E\to F$ compact. 
%\end{Teor}
The first major result for the real method appeared when Hayakawa \cite{Hayakawa1969} proved the following.
\begin{Teor}[Hayakawa]\label{hayakawa}
    Let $(E_0,E_1),(F_0,F_1)$ be two compatible couples of Banach spaces and $T:(E_0,E_1)\to (F_0,F_1)$ an admissible operator such that $T:E_i\to F_i$ is compact, $i=0,1$. Then, for any $1\leq q<\infty$, $0<\theta<1$, $T:(E_0,E_1)_{\theta,q}\to (F_0,F_1)_{\theta,q}$ compactly.
\end{Teor}
This result says that the real method for $1\leq q<\infty$ preserves two-sided compactness without any additional hypotheses on the couples. However, the case $q=\infty$ was still missed. Moreover, the arguments of Hayakawa were very difficult to follow and many details were unclear. Two decades later, in \cite{CobosEdmundsPotter1990} a transparent proof of the results of Hayakawa, extending it to the case $q=\infty$, and also to the quasi-Banach case $0<q<1$, was given.

\begin{Teor}[Cobos-Edmunds-Potter]
    Let $(E_0,E_1),(F_0,F_1)$ be two compatible couples of Banach spaces and $T:(E_0,E_1)\to (F_0,F_1)$ an admissible operator such that $T:E_i\to F_i$ is compact, $i=0,1$. Then, for any $0< q\leq\infty$, $0<\theta<1$, $T:(E_0,E_1)_{\theta,q}\to (F_0,F_1)_{\theta,q}$ compact.
\end{Teor}
Their ideas also generalized the results of Lions and Peetre in the following way.
\begin{Teor}[Cobos-Edmunds-Fern\'andez-Potter] Let $(E_0,E_1), (F_0,F_1) $ be two compatible couples of Banach spaces, $0<\theta<1$ and $1\leq q\leq \infty$. Assume that $T:(E_0,E_1)\to (F_0,F_1)$ is an addmisible operator.
\begin{itemize}
    \item If $T:E_0\to F_0$ compact and $F_1\xhookrightarrow{}F_0$, then $T:(E_0,E_1)_{\theta,q}\to (F_0,F_1)_{\theta,q}$ compact.
    \item If $T:E_1\to F_1$ compact and $E_1\xhookrightarrow{}E_0$, then $T:(E_0,E_1)_{\theta,q}\to (F_0,F_1)_{\theta,q}$ compact.
\end{itemize}
\end{Teor}
The second item was proved by Cobos and Fern\'andez in \cite{CobosFernandez1989}. Finally, in 1992, the one-sided question for the real method was solved by Michael Cwikel \cite{Cwikel1992}.
\begin{Teor}[Cwikel]\label{Cwikel}
    Let $(E_0,E_1),(F_0,F_1)$ be two compatible couples of Banach spaces and $T:(E_0,E_1)\to (F_0,F_1)$ an admissible operator such that $T:E_0\to F_0$ compact. Then, for any $1\leq q\leq\infty$, $0<\theta<1$, $T:(E_0,E_1)_{\theta,q}\to (F_0,F_1)_{\theta,q}$ is compact.
\end{Teor}
Soon, the result of Cwikel was generalized to the abstract setting of Aronszajn-Gagliardo functors (see \cite{AronszajnGagliardo1965}) by Cobos, K\"uhn and Schonbek \cite{CobosKuhnSchonbek1992}. More generalizations of the compactness and real method has been established since then, like the one-sided result for the quasi-Banach setting \cite{CobosPersson1998} motivated by the Krasnoselskii's theorem for the quasi-Banach case \cite{ZabreikoPustylnik1965}. Very important results have been obtained for the compactness on the limiting cases $\theta=0,1$, where variations of the real method like the logarithmic one were introduced and developed  in \cite{CobosSegurado2015},\cite{CobosLuzKunUllrich2009},\cite{FernandezMartinezSeguradoSignes2016}. We refer to \cite{GarciaSaez2024} and \cite{Cobos2009} to a larger exposition on the classical problem and its development.\\

The compactness question for the complex interpolation method has been the subject of extensive research, but its development has followed a distinct path. In fact, the one-sided and even the two-sided results are still open problems. In the seminal paper of Calder\'on \cite{Calderon1964}, it is implicit that compactness results hold for complex interpolation under some assumption on the couple. Partial results were given, for example as a consequence of functorial Lions-Peetre theorems, since they apply to interpolation functors of exponent $\theta$,  but significant advances in the problem were absent for several years. Interest in the problem was rekindled in 1989 with the publication of \cite{CobosPeetre1989}, which introduced powerful new ideas. These ideas were combined with the results of Janson \cite{Janson1981} by Cwikel in \cite{Cwikel1992} to prove that the problem was, essentially, equivalent to show if the compactness is preserved for a particular choice of the couples (some vector-valued sequence spaces). The next important result was obtained by Cobos, K\"uhn and Schonbek in \cite{CobosKuhnSchonbek1992} for Banach lattices of complex-valued functions (we refer to section 4 in \cite{RaynaudTradacete2010} to a brief discussion on the topic and more references).

% Moreover, suppose that there exists a directed set of linear operators $\{T_\lambda\}$ in $F_0+F_1$ such that $T_\lambda(E_0)\subset E_0$ with the norms of $T_\lambda:F_j\to F_j$ uniformly bounded, $j=0,1$; the spaces $T_\lambda(E_0)$ are finite-dimensional and $T_\lambda x\to x$ for every $x\in E_0$ in the convergence sense of nets. In this case, $T:[E_0,E_1]_\theta\to [F_0,F_1]_\theta$ compactly. 
\begin{Teor}[Cobos-K\"uhn-Schonbek]
    Let $(E_0,E_1)$ and $(F_0,F_1)$ be compatible couples of Banach spaces formed by complexified Banach lattices of measurable functions defined on $\sigma$-finite measure spaces. Assume that $T:(E_0,E_1)\to (F_0,F_1)$ is an admissible operator such that $T:E_0\to F_0$ is compact. Then, if $F_0$ or $F_1$ has absolutely continuous norm or has the Fatoy propery, then $T:[E_0,E_1]_\theta\to [F_0,F_1]_\theta$ is compact for $\theta\in (0,1)$.
\end{Teor}
In the same work, an interesting extrapolation result was given:
\begin{Teor}
    Let $(E_0,E_1),(F_0,F_1)$ be compatible couples of Banach spaces and $T:(E_0,E_1)\to (F_0,F_1)$ an admissible operator. Suppose that for some $\theta_0\in (0,1)$ $T:[E_0,E_1]_{\theta_0}\to [E_0,E_1]_{\theta_0}$ is compact. Then, $T:[E_0,E_1]_{\theta}\to [E_0,E_1]_{\theta}$ is compact for every $\theta\in (0,1)$.
\end{Teor}
The major result on the compactness problem for the complex method was given by Cwikel and Kalton in \cite{CwikelKalton1995} and lately extended by Cwikel \cite{cwikel2008}:
\begin{Teor}[Cwikel-Kalton]\label{ComplexCompact}
    Let $(E_0,E_1),(F_0,F_1)$ compatible couples of Banach spaces, $T:(E_0,E_1)\to (F_0,F_1)$ an admissible operator such that $T:E_0\to F_0$ is compact, $\theta\in (0,1)$ and $E,F$ two Banach spaces. Then, if one of the following cases holds:
    \begin{itemize}
        \item $E_0$ has the UMD property.
        \item $E_0$ is reflexive and $E_0=[E,E_1]_{\theta_0}$, $\theta_0\in (0,1)$.
        \item $F_0=[F,F_1]_{\theta_1}, \theta_1\in (0,1)$.
        \item $E_0$ and $E_1$ are complexified Banach lattices of measurable functions over a common measure space
        \item $F_0$ and $F_1$ are complexified Banach lattices of measurable functions over a common measure space such that at least one of them has either absolutely continuous norm or the Fatou property.
    \end{itemize} Then, $T:[E_0,E_1]_\theta\to [F_0,F_1]_\theta$ is compact.
\end{Teor}
Several equivalent forms of the problem were established by Cwikel, Krugljak and Mastylo in \cite{CwikelKrugljakMastylo1996}.
Since this major result, many interesting partial cases have been solved, for example see the works of Schonbek \cite{Schonbek2000} and Cobos, Fern\'andez-Cabrera and Mart\'inez \cite{CobosFernandezCabreraMartinez2004}. For a discussion on the problem and some recent advances we refer to \cite{CwikelRochberg2014}. \\

We recall that we say that a Banach space $E$ satisfies the UMD condition (unconditional martingale differences) if 
and only if the Hilbert transform extends to a bounded operator on $L^p(
\R^n; E)$. For a complete introduction to those spaces we refer to \cite[Chapter 4 \& 5]{UMD2016}. In the following section, we will prove compact embeddings using Theorem \ref{ComplexCompact} and the fact that Sobolev spaces are UMD (see \cite[Example 4.2.18]{UMD2016}).
\subsection{Bessel potential spaces}

As mentioned in the introduction, Bessel potential spaces are classical fractional Sobolev spaces of great importance in the study of partial differential equations. Classical definition of Bessel potential spaces is the following.

\begin{Defi}\em[Bessel potential space]
    Let $s\in \R$ and $1\leq p\leq \infty$. We define the \textit{Bessel potential space} $\Lambda^{s,p}(\R^n)$ as $$\Lambda^{s,p}(\R^n):=\{f \in \mathcal{S}(\R^n): \Lambda_{-s}f\in L^p(\R^n)\},$$ with the norm $$\norm{f}_{\Lambda^{s,p}}:=\norm{\Lambda_{-s}f}_p,$$ 
    where $\Lambda_s$ is the Bessel potential of order $s$ of $f$ defined as $$\Lambda_{s}f=\mathfrak{F}^{-1}\left((1+4\pi^2|\xi|^2)^{-s/2}\mathfrak{F}f(\xi)\right),$$ where $\mathfrak{F}$ denotes the Fourier transform.
\end{Defi}

For functions in the Schwartz class $\mathcal{S}(\R^n)$, Bessel potential is clearly well defined. We could also define complex interpolation of Sobolev spaces in the following way. 
\begin{Defi}\em[Complex interpolation of Sobolev spaces]
    Let $\theta\in [0,1]$, $1<p<\infty$ and $k\in \mathbb{Z}$. Let $s=k\theta$. We define the space $H^{s,p}(\R^n)$ as the complex interpolation space $$H^{s,p}(\R^n):=[L^p(\R^n),W^{k,p}(\R^n)]_\theta,$$ with the identification $H^{0,p}(\R^n)=L^p(\R^n)$. In particular, for $s\in [0,1]$, $$H^{s,p}(\R^n)=[L^p(\R^n),W^{1,p}(\R^n)]_s.$$
\end{Defi}
The crucial mathematical fact, known in the mathematical analysis community since the sixties of the past century, is that for every $s\in\R$ and $1<p<\infty$, $$\Lambda^{s,p}(\R^n)=H^{s,p}(\R^n),$$ with equivalence of the norms (a proof of this fact and references to the classical literature can be found in \cite{BellidoGarciaA2025}). Using this equivalence and the theory of complex interpolation, many interesting properties of Bessel potential spaces can be obtained. In \cite{BellidoGarciaA2025}, this connection of Bessel potential spaces has been revisited, providing direct proofs based on interpolation theory for all the results collected in the following theorem. 
\begin{Teor}\label{BesselProps}
    Let $s\in (0,1)$ and $1<p<\infty$. Then,
    \begin{itemize}
        \item $H^{s,p}(\R^n)$ is a complete, separable and reflexive normed space.
        \item $C_c^\infty(\R^n)$ and $W^{1,p}(\R^n)$ are dense in $H^{s,p}(\R^n)$.
        \item For $0<s_0<s_1<1$, $H^{s_1,p}(\R^n)\xhookrightarrow{}H^{s_0,p}(\R^n)$.
        \item The following continuous embeddings hold 
        $$\begin{cases}
            H^{s,p}(\R^n)\xhookrightarrow{}L^{p_s^*}(\R^n),\,p_s^*=\frac{np}{n-sp},&sp<n,\\
            H^{s,p}(\R^n)\xhookrightarrow{} \normalfont{\textbf{BMO}}(\R^n),L^q(\R^n),\,p\leq q<\infty&sp=n,\\
            H^{s,p}(\R^n)\xhookrightarrow{}C^{0,\mu_s}(\R^n),\,\mu_s=s-n/p,&sp>n.
        \end{cases}$$
        \item Let $s_0,s_1\in \R^n$ and $1<p_0<p_1<\infty$, then $$[H^{s_0,p_0}(\R^n),H^{s_1,p_1}(\R^n)]_{\theta}=H^{s(\theta),p(\theta)}(\R^n),$$ where $s(\theta)=(1-\theta)s_0+\theta s_1,\,\frac{1}{p(\theta)}=\frac{1-\theta}{p_0}+\frac{\theta}{p_1}.$
        \item Let $0<t<s<1$, $1<p<\infty$ and $1<p<q\leq \frac{np}{n-(s-t)p}$. Then, $$H^{s,p}(\R^n)\xhookrightarrow{}H^{t,q}(\R^n).$$
    \end{itemize}
\end{Teor}

Another scale of fractional Sobolev spaces is the Sobolev-Slobodeckij, or Gagliardo, spaces. These spaces are extensively utilized due to their strong connection with the fractional $p$-Laplacian and are often termed as fractional Sobolev spaces. To avoid confusion, since Bessel potential spaces are also fractional-order Sobolev spaces, we will refer to them specifically as Gagliardo spaces.. 
\begin{Defi}[Gagliardo spaces]\em
    Let $s\in (0,1)$ and $p\in [1,\infty)$. We define the Gagliardo space $W^{s,p}(\R^n)$ as $$W^{s,p}(\R^n):=\Bigg\{u\in L^p(\R^n): \frac{|u(x)-u(y)|}{|x-y|^{n+sp}}\in L^p(\R^n\times \R^n)\Bigg\},$$ with the norm $$\norm{u}_{W^{s,p}}:=\norm{u}_p+[u]_{W^{s,p}},$$ where $$[u]_{W^{s,p}}:=\int_{\R^n}\int_{\R^n}\frac{|u(x)-u(y)|^p}{|x-y|^{n+sp}}\,dx\,dy,$$ is the Gagliardo seminorm. In the case $p=\infty$, we have the natural identification of  $W^{s,\infty}$ with $C^{0,s}$.
\end{Defi} 
As well as Bessel potential spaces, Gagliardo spaces are also obtained by means of interpolation methods, but in this case, real interpolation of Sobolev spaces. In particular, $$W^{s,p}(\R^n)=\left(L^p(\R^n),W^{1,p}(\R^n)\right)_{s,p},\,s\in (0,1),\,1<p<\infty,$$ with equivalence of norms. We refer to \cite[Theorem~12.5]{Leoni2023} and \cite[Example~1.8]{Lunardi2018} for detailed proofs of this fact, and to \cite{DiNezzaPalatucciValdinoci2012} for a complete discussion of Gagliardo spaces without appealing to interpolation theory. 

%With analogous interpolation arguments as the ones for Bessel potential spaces in \cite{BellidoGarcia2024}, the following properties of Gagliardo spaces are obtained:
%\begin{Teor}\label{GagliardoProps}
   % Let $s\in (0,1)$ and $1<p<\infty$. Then,
   % \begin{itemize}
   %     \item $W^{s,p}(\R^n)$ is a complete, separable and reflexive normed space.
   %     \item $C_c^\infty(\R^n)$ and $W^{1,p}(\R^n)$ are dense in $W^{s,p}(\R^n)$.
   %     \item For $0<s_0<s_1<1$, $W^{s_1,p}(\R^n)\xhookrightarrow{}W^{s_0,p}(\R^n)$.
   %     \item $$\begin{cases}
   %         W^{s,p}(\R^n)\xhookrightarrow{}L^{p_s^*}(\R^n),\,p_s^*=\frac{np}{n-sp},&sp<n,\\
   %         W^{s,p}(\R^n)\xhookrightarrow{} \textbf{BMO}(\R^n),L^q(\R^n),\,p\leq q<\infty&sp=n,\\
   %         W^{s,p}(\R^n)\xhookrightarrow{}C^{0,\mu_s}(\R^n),\,\mu_s=s-n/p,&sp>n.
   %     \end{cases}$$
   %     \item Let $0<s_0<s_1\leq1$ and $1<p_0<p_1<\infty$, then $$\left(W^{s_0,p_0}(\R^n),W^{s_1,p_1}(\R^n)\right)_{\theta,p(\theta)}=[W^{s_0,p_0}(\R^n),W^{s_1,p_1}(\R^n)]_{\theta}=W^{s(\theta),p(\theta)}(\R^n),$$ where $s(\theta)=(1-\theta)s_0+\theta s_1,\,\frac{1}{p(\theta)}=\frac{1-\theta}{p_0}+\frac{\theta}{p_1}.$ Also, $$\left(L^{p_0}(\R^n),W^{s,p_1}(\R^n)\right)_{\theta,p(\theta)}=[L^{p_0}(\R^n),W^{s,p_1}]_\theta=W^{\theta s,p(\theta)}(\R^n).$$
    %    \item Let $0<t<s<1$, $1<p<\infty$ and $1<p<q\leq \frac{np}{n-(s-t)p}$. Then, $$W^{s,p}(\R^n)\xhookrightarrow{}W^{t,q}(\R^n).$$
   % \end{itemize}
%\end{Teor}

This identification with real interpolation spaces reveals a strong relationship between Gagliardo and Bessel spaces, despite their apparent differences. By applying general interpolation arguments, the following embeddings that connect Bessel and Gagliardo spaces are established. For detailed proofs, we refer interested readers to \cite{BellidoGarciaA2025}.
\begin{Teor}[Bessel versus Gagliardo]\label{Contiguity}
    Let $0<s_0<s<s_1<1$ and $1<p<\infty$. Then, \begin{itemize}
        \item $W^{s,2}(\R^n)=H^{s,2}(\R^n)$ with equivalence of the norms.
        \item $H^{s_1,p}(\R^n)\xhookrightarrow{}W^{s,p}(\R^n)\xhookrightarrow{}H^{s_0,p}(\R^n).$
        \item $$\begin{cases}
        W^{s,p}(\R^n)\xhookrightarrow{}H^{s,p}(\R^n),&1<p\leq 2,\\
        H^{s,p}(\R^n)\xhookrightarrow{}W^{s,p}(\R^n),&2\leq p<\infty,
    \end{cases}$$ with strict inclusions unless $p=2$.
    \end{itemize}
\end{Teor}

Since the aim of this paper is to obtain compact embeddings of Bessel potential spaces, we have to work on bounded domains, and we need to define in a precise way Bessel potential spaces on domains. Let $s\in (0,1)$, $1<p<\infty$ and $\Omega\subset \R^n$ be an open set. Let $f$ a function defined on $\R^n$. We define the restriction operator to $\Omega$, $R_\Omega$, as the operator $R_\Omega f=f|_\Omega$. We define the Bessel potential space on $\Omega$ in the natural manner, that is to say, as a quotient space,
$$\Lambda^{s,p}(\Omega):=R_\Omega\left(\Lambda^{s,p}(\R^n)\right)=\Lambda^{s,p}(\R^n)/\sim,$$ 
where, $f\sim g$ if and only if $R_\Omega f(x)=R_\Omega g(x)$ a.e. in $\Omega$. The norm of the space is the quotient norm, 
$$\norm{f}_{\Lambda^{s,p}(\Omega)}=\operatorname{inf}\{\norm{g}_{\Lambda^{s,p}}: R_\Omega g=f\}.$$ As in the case of $\Omega=\R^n$, we want to have the equivalence of such spaces with the complex interpolation scale, that may be defined also for Sobolev spaces on domains, 
$$H^{s,p}(\Omega):=[L^p(\Omega),W^{1,p}(\Omega)]_s,\,s\in (0,1),\,1<p<\infty$$ equipped with the interpolation norm. For suitable domains $\Omega$ one would expect to have $$H^{s,p}(\Omega)=\Lambda^{s,p}(\Omega),$$ with equivalence of the norms. It is known (see \cite[Proposition 2.4]{jerison1995}) that when $\Omega$ is a bounded set with a Lipschitz boundary, the spaces $H^{s,p}(\Omega)$ form an interpolation scale. To ensure completeness in this manuscript, we include a proof of the equality between $H^{s,p}(\Omega)$ and $\Lambda^{s,p}(\Omega)$. The key is the existence of an extension operator for such $\Omega$ \cite[Theorem VI.5]{Stein1970}. We utilize the identification $W^{0,p}(\Omega) = L^p(\Omega)$.

\begin{Teor}[Extension operator]\label{SteinExtension}
    Let $k\ge 0$ be an integer, $1\le p \le +\infty$, and $\Omega\subset \R^n$ a Lipschitz open set. Then, there exists an extension operator, i.e., a bounded linear operator $\mathcal{S}\in \mathcal{L}(W^{k,p}(\Omega),W^{k,p}(\R^n))$ such that  $\left(\mathcal{S}u\right)|_{\Omega}=u$, for every $u\in W^{k,p}(\Omega)$.
\end{Teor}
The crucial fact about the extension operator for our aim here is that it does not depend on $k$ for the spaces $W^{k,p}$ (it does not depend on $p$ either).  
\begin{Teor}\label{BesselDomain}
    Let $s\in (0,1)$ and $1<p<\infty$. Then, $\mathcal{S}$ is an extension operator from $H^{s,p}(\Omega)$ to $H^{s,p}(\R^n)$, i.e., $\mathcal{S}\in\mathcal{L}(H^{s,p}(\Omega),H^{s,p}(\R^n))$ and $\left(\mathcal{S}u\right)|_{\Omega}=u$, for every $u\in H^{s,p}(\Omega)$. Moreover, 
    $$H^{s,p}(\Omega)=\Lambda^{s,p}(\Omega),$$ 
    and the norms of both spaces are equivalent. 
\end{Teor}
\noindent{}\textbf{Proof:} The restriction operator $$R_\Omega:L^p(\R^n)+W^{1,p}(\R^n)\to L^p(\Omega)+W^{1,p}(\Omega),$$ is an admissible operator with norm less or equal than $1$, and by complex interpolation $$R_\Omega:H^{s,p}(\R^n)\to H^{s,p}(\Omega),$$ is a bounded linear operator with norm less or equal than $1$. The extension operator of {Theorem \ref{SteinExtension}}, $$\mathcal{S}:L^p(\Omega)+W^{1,p}(\Omega)\to L^p(\R^n)+W^{1,p}(\R^n),$$ is an admissible operator and by complex interpolation $$\mathcal{S}:H^{s,p}(\Omega)\to H^{s,p}(\R^n),$$ is an extension operator. It follows now that \begin{align*}
    \mathcal{S}\left(H^{s,p}(\Omega)\right)&\subset H^{s,p}(\R^n)=\Lambda^{s,p}(\R^n),\\
R_\Omega\left(\Lambda^{s,p}(\R^n)\right)&=R_\Omega\left(H^{s,p}(\R^n)\right) \subset H^{s,p}(\Omega),
\end{align*} and since $R_\Omega\circ \mathcal{S}=\operatorname{id}$, we conclude that 
$$H^{s,p}(\Omega)\subset R_\Omega\left(\Lambda^{s,p}(\R^n)\right)\subset H^{s,p}(\Omega),$$ and hence $$H^{s,p}(\Omega)=\Lambda^{s,p}(\Omega).$$ Now we will prove the equivalence of the norms. Set 
$$C_\mathcal{C}:=\| \mathcal{S}\|_{\mathcal{L}(H^{s,p}(\Omega),H^{s,p}(\R^n))},$$
and $C_1,C_2>0$ such that
\[C_2 \|\cdot\|_{H^{s,p}}\le \|\cdot\|_{\Lambda^{s,p}}\le C_1 \|\cdot\|_{H^{s,p}}.\]
Let $u\in H^{s,p}(\Omega)$. Then, $$\norm{u}_{\Lambda^{s,p}(\Omega)}\leq \norm{\mathcal{S}u}_{\Lambda^{s,p}}\leq C_1\norm{\mathcal{S}u}_{H^{s,p}}\leq C_1C_\mathcal{S}\norm{u}_{H^{s,p}(\Omega)}.$$ Now, let $u\in \Lambda^{s,p}(\Omega)$. There exists $v\in \Lambda^{s,p}(\R^n)$ such that $v|_\Omega=u$, and hence $$\norm{u}_{H^{s,p}(\Omega)}= \norm{R_\Omega v}_{H^{s,p}}\leq \norm{v}_{H^{s,p}}\leq C_2^{-1}\norm{v}_{\Lambda^{s,p}}.$$ Now, taking the infimum over all possible functions $v$ it yields that $$\norm{u}_{H^{s,p}(\Omega)}\leq C_2^{-1}\norm{u}_{\Lambda^{s,p}(\Omega)}.\qed$$

Previous definitions and results are easily extended to the case $s>0$, and hence we have the definition of $H^{s,p}(\Omega)$ as an interpolation space and its equivalence with $\Lambda^{s,p}(\Omega)$ for $s>0$ when $\Omega$ is a Lipschitz domain. 

The assumption on the domain $\Omega$ to be Lipschitz can be relaxed to the more general requirement of being a locally uniform domain. A domain $\Omega\subset \R^n$ is $(\varepsilon,\delta)$-locally uniform if between any pair of points $x,y\in \Omega$ such that $|x-y|<\delta$, there is a rectifiable arc $\gamma\subset\Omega$ of length at most $|x-y|/\varepsilon$ and having the property that for all $z\in \gamma$, 
$$\mbox{dist}(z,\partial \Omega)\ge \frac{\varepsilon|z-x||z-y|}{|x-y|}.$$ 
This family of domains contains uniformly Lipschitz domains as special cases. In \cite{ROGERS2006619}, it is shown that for every $(\varepsilon,\delta)$-locally uniform domain $\Omega$, there exists a linear bounded operator  $u\mapsto \mathcal{R}u$ such that for any $k\in\mathbb{N}\cup\{0\}$ and $1\leq p\leq \infty$, $\mathcal{R}$ is an extension operator $W^{k,p}(\Omega)\to W^{k,p}(\R^n)$ with norm depending on $n,k,p,\varepsilon,\delta$ Then, by analogous interpolation arguments as in Theorem \ref{BesselDomain}, we would get the equivalence 
$$H^{s,p}(\Omega)=\Lambda^{s,p}(\Omega),$$ 
for such a more general family of domains.

\subsection{Riesz fractional gradient}

As mentioned in the introduction, possibly the most important reason for this renewed interest in Bessel potential spaces is their connection with the Riesz fractional gradient, since the pioneering works \cite{ShiehSpector2015,ShiehSpector2018}. The Riesz fractional gradient of order $s\in(0,1)$ for a function $u\in C_c^\infty(\R^n)$ is defined as
\[D^s u(x)=c_{n,s} \int_{\mathbb{R}^n} \frac{u(x)-u(y)}{|x-y|^{n+s}}\frac{x-y}{|x-y|}\,dy,\]
where 
$$c_{n,s}=(n+s-1)\frac{\Gamma\left(\frac{n+s-1}{2}\right)}{\pi^\frac{n}{2} 2^{1-s}\Gamma\left(\frac{1-s}{2}\right)}, $$
a normalization constant, with $\Gamma$ the Euler gamma function. In \cite{ShiehSpector2015,ShiehSpector2018}, it was demonstrated that Bessel potential spaces coincide, with equivalence of norms, with the closure of $C_c^\infty(\mathbb R^n)$ with respect to the norm
$$\norm{u}_{D^s,p}=\|u\|_{L^p(\mathbb{R}^n)}+\|D^s u\|_{L^p(\mathbb{R}^n,\mathbb{R}^n)}.$$
Moreover, the Riesz fractional gradient exhibits several notable properties. In \cite{Silhavy2020}, the uniqueness of the Riesz fractional gradient, up to a multiplicative constant, was established under natural conditions (invariance under translations and rotations, homogeneity under dilatation, and certain continuity in the sense of distributions). Furthermore, for a Sobolev function $u \in W^{1,p}(\mathbb{R}^n)$, $D^s u$ converges strongly to $Du$ in $L^p(\mathbb R^n)$ \cite{BellidoCuetoMoraCorral2021}. Additionally, defining the fractional divergence of order $s$ of a function $\psi:\R^n\to \R^n$ as
\[\operatorname{div}^s \psi(x)=-c_{n,s}\int_{\R^n} \frac{\psi(x)+\psi(y)}{|x-y|^{n+s}}\cdot \frac{x-y}{|x-y|}\,dy, \]
integration by parts holds (see \cite{BellidoCuetoMoraCorral2020} and references therein). Moreover, fractional Laplacian in $\mathbb{R}^n$ can be obtained as the fractional divergence of the Riesz fractional gradient. Consequently, the Riesz fractional gradient $D^s u$ is a truly fractional differential object, and this, combined with the identification of Bessel potential spaces in \cite{ShiehSpector2015}, makes Bessel potential spaces an intriguing and valuable framework for the analysis of fractional PDEs.

In \cite{ShiehSpector2015}, beyond the connection of Bessel potential spaces with Riesz fractional gradients mentioned above, it was shown a representation theorem called fractional fundamental theorem of calculus, which establishes that a function can be recovered from its Riesz fractional gradient through convolution with a  certain kernel. That is a remarkable and important result, and, as a consequence, continuous embeddings of Bessel potential spaces were given in \cite{ShiehSpector2015} and compactness of those embeddings in \cite{ShiehSpector2018}. Those results are based on the identification of Bessel potential spaces as the spaces of $L^p$ functions whose Riesz fractional gradient is also a $L^p$ function. In this paper, we provide a proof of compactness for continuous embeddings given in \cite{BellidoGarciaA2025} based on estimations on translations of given smooth compactly supported functions based on the Riesz fractional gradient. For readers' convenience, we include here a new simple proof of the fractional fundamental theorem of calculus. 

\begin{Prop}\label{FTOC}
    Let $s\in (0,1)$ and $u\in C_c^\infty(\R^n)$. Then, $$u(x)=c_{n,-s}\int_{\R^n}D^su(y)\cdot \frac{x-y}{|x-y|^{n-s+1}}\,dy.$$
\end{Prop}
\noindent{}\textbf{Proof:}
Let us define 
$$K_s(x):=c_{n,-s}\frac{x}{|x|^{n-s+1}},$$ 
for $x\in\R^n\setminus\{0\}$. We claim that 
$$\widehat{K}_s(\xi)=-i\frac{\xi}{|\xi|}|2\pi\xi|^{-s}.$$ 
Then, for $u\in C_c^\infty(\R^n)$, 
$$c_{n,-s}\int_{\R^n}D^su(y)\cdot \frac{x-y}{|x-y|^{n-s+1}}\,dy=D^su*K_s(x),$$ so 
$$\widehat{D^su*K(x)}(\xi)=\widehat{D^su}\cdot \widehat{K_s}(\xi)=\frac{2\pi i\xi}{|2\pi \xi|^{1-s}}\widehat{u}(\xi)\cdot (-i)\frac{\xi}{|\xi|}|2\pi\xi|^{-s}=\widehat{u}(\xi),$$ 
hence 
$$u(x)=c_{n,-s}\int_{\R^n}D^su(y)\cdot \frac{x-y}{|x-y|^{n-s+1}}\,dy.$$ 
We only have to compute the Fourier transform of $K_s$. In order to do this, defining
\[\gamma_{1+s}=\pi^{\frac{n}{2}} 2^{1+s} \frac{\Gamma\left(\frac{1+s}{2}\right)}{\Gamma\left(\frac{n-1-s}{2}\right)},\]
observe that 
$$D\left(\frac{1}{\gamma_{1+s}}\frac{1}{|x|^{n-s-1}}\right)=-\frac{(n-s-1)}{\gamma_{1+s}}\frac{x}{|x|^{n-s+1}}=-c_{n,-s}\frac{x}{|x|^{n-s+1}}=-K_s(x),$$ 
so
$$\mathfrak{F}\left\{D\left(\frac{1}{\gamma_{1+s}}\frac{1}{|x|^{n-s-1}}\right)\right\}(\xi)=-\widehat{K_s}(\xi).$$ Since $$\widehat{DI_{1+s}}(\xi)=2\pi i\xi \widehat{I_{1+s}}(\xi)=2\pi i\xi |2\pi\xi|^{-(1+s)}=\frac{i\xi}{|\xi|}|2\pi\xi|^{-s},$$ hence $$\widehat{K_s}(\xi)=-\frac{i\xi}{|\xi|}|2\pi\xi|^{-s},$$ as we wanted to prove. \qed\\

\subsection{Three classical results}

For the sake of completeness in the exposition,  we recall three very classical results of functional analysis that we use in the following section. 

\begin{Teor}[Ascoli-Arzel\'a]\label{ascoliarzela}
    Let $(X,d)$ a compact metric space and $W$ a uniformly bounded subset of $\mathcal{C}(X)$, the space of continuous functions from $X$ to $\R$. If $W$ is equicontinuous, i.e., for every $\varepsilon>0$ there exists $\delta>0$ such that if $d(x,y)<\delta$, then
    $$|f(x)-f(y)|<\varepsilon,$$
    for all, then $W$ is relatively compact in $C(X)$.
\end{Teor}

\begin{Teor}[Fr\'echet-Kolmogorov-Riesz]\label{FrechetKolmogorov}
    Let $W$ a bounded set in $L^p(\R^n)$ for $1\leq p<\infty$ such that $$\lim_{|h|\to 0^+}\norm{u(\cdot+h)-u(\cdot)}_p=0,$$ uniformly in $u\in W$, i.e., for every $\varepsilon>0$, there exists $\delta>0$ such that 
    $$\norm{u(\cdot+h)-u(\cdot)}_p<\varepsilon,\,\,|h|<\delta,$$
    for all $u\in W$. Then, for any subset $\Omega\subset \R^n$ with finite measure, $W_{|\Omega}=\left\{ f_{|\Omega}\;:\;f\in W\right\} $ is relatively compact in $L^p(\Omega)$. 
\end{Teor}
See \cite[Theorems~4.25 and 4.26]{Brezis2010} for the proofs of the theorems. The following proposition is elementary. 

\begin{Prop}
    Let $E,F,G$ three Banach spaces and $T\in \mathcal{L}(E,F)$ and $S\in\mathcal{L}(F,G)$ such that at least one of them is compact (the image of the operator is compact) Then, the composition $S\circ T:E\to G$ is compact.
\end{Prop}

%\noindent{}\textbf{Proof:} Suppose that $T:E\to F$ compactly and $S:F\to G$ continuously (the other case is analogous). Let $B_E$ the closed unit ball in $E$. Since $T$ %is compact, $T(B_E)$ is relatively compact, i.e., $\overline{T(B_E)}$ is a compact subset of $F$. Since continuous images of compact sets are compact, %$S\left(\overline{T(B_E)}\right)$ is compact, and since $ST(B_E)\subseteq S\left(\overline{T(B_E)}\right)$, we conclude that $ST:E\to G$ compactly.\qed

Given two Banach spaces $E$ and $F$, the notation $E\dhookrightarrow F$ means that $E$ is compactly embedded into $F$, i.e., bounded subsets of $E$ are relatively compact with respect to the topology in $F$.  
\begin{Coro}\label{compactcompos}
    Let $E,F,G$ three Banach spaces such that $E\dhookrightarrow F$ and $F\xhookrightarrow{}G$ (respectively, $E\xhookrightarrow{}F$ and $F\dhookrightarrow G$). Then, $E\dhookrightarrow G$.
\end{Coro}

\section{Compact embeddings}

In this section we give three different proofs of Theorem \ref{BesselCompact}. The first one, for the case $1< p< \infty$, is based on known compactness results on the complex method of interpolation, in particular Theorem \ref{ComplexCompact}. The second one, for the case $1\le p <\infty$, is based on estimations on translations involving the Riesz fractional gradient and applying the Fréchet-Kolmogorov theorem. The third one, just uses the contiguity between Gagliardo and Bessel potential spaces, which also follows from interpolation techniques, and the compactness results for the former. 

\subsection{A proof based on interpolation theory}

In this section we prove Theorem \ref{BesselCompact} for $1<p<\infty$ and another interesting result on compact embeddings between Bessel spaces. We begin with the proof of the theorem.

%\begin{Teor}[Fractional Rellich-Kondrachov-Subcritical case]\label{Weak subcritical}
 %   Let $\Omega$ a subset of $\R^n$ with finite measure. Let $s\in (0,1)$ and $1<p<\infty$ such that $sp<n$. Then, $$H^{s,p}%(\R^n)\xhookrightarrow{}\xhookrightarrow{}L^{q}(\Omega),\,1<q<p_s^*=\frac{np}{n-sp}.$$
%\end{Teor}
%\noindent{}\textbf{Proof:} 

\begin{itemize}
\item[a)] We first suppose that $p<n$, and hence $p^*=\frac{np}{n-p}$ is well defined. By the Rellich-Kondrachov theorem for classical Sobolev spaces, $$W^{1,p}(\R^n)\dhookrightarrow L^q(\Omega),$$ for every $1\leq q<p^*$. Let $q\in [p,p^*)$ arbitrary and consider the inclusion operator $$i:L^p(\R^n)+W^{1,p}(\R^n)\to L^p(\Omega)+L^q(\Omega),$$ which is an admissible operator such that $i:W^{1,p}(\R^n)\to L^q(\Omega)$ is compact. Then, since $W^{1,p}(\R^n)$ is a UMD-space, $$i:[L^p(\R^n),W^{1,p}(\R^n)]_s\to [L^p(\Omega),L^q(\Omega)]_s,$$ is a compact inclusion by {Theorem \ref{ComplexCompact}}. We have that $[L^p(\R^n),W^{1,p}(\R^n)]_s=H^{s,p}(\R^n)$ and $[L^p(\Omega),L^q(\Omega)]_s=L^{r_s(q)}(\Omega),$ where $$\frac{1}{r_s(q)}=\frac{1-s}{p}+\frac{s}{q}.$$ The function $q\mapsto r_s(q)$ is a continuous function such that $r_s(p)=p$ and $$\lim_{q\to \left(p^*\right)^-}r_s(q)=\left(\frac{1-s}{p}+\frac{s}{p^*}\right)^{-1}=\left(\frac{1-s}{p}+\frac{s(n-p)}{np}\right)^{-1}=p_s^*=\frac{np}{n-sp},$$ 
and hence, it maps continuously $[p,p^*)$ into $[p,p_s^*)$, so $$H^{s,p}(\R^n)\dhookrightarrow L^r(\Omega),\,p\leq r<p_s^*.$$ 
Now, since $\Omega$ is bounded, $L^{q_2}(\Omega)\xhookrightarrow{}L^{q_1}(\Omega)$ for $q_1<q_2$, thus, for $q\in [1,p)$, 
$$H^{s,p}(\R^n)\dhookrightarrow{} L^p(\Omega)\xhookrightarrow{}L^q(\Omega),$$ 
and by {Corollary \ref{compactcompos}}, 
$$H^{s,p}(\R^n)\dhookrightarrow{}L^q(\Omega),$$ 
for all $1\leq q<p_s^*$. Now, suppose that $p=n$, since $W^{1,p}(\R^n)\dhookrightarrow L^q(\Omega)$ for $q\in [1,+\infty)$, the same argument as above shows that 
%$$H^{s,p}(\R^n)\dxhookrightarrow L^{r_s(q)}(\Omega),$$ for $q\in [p,+\infty)$. Again, $r_s(p)=p$ and 
$$\lim_{q\to +\infty}r_s(q)=\lim_{q\to +\infty}\left(\frac{1-s}{p}+\frac{s}{q}\right)^{-1}=\frac{p}{1-s}=\frac{np}{n-ns}=\frac{np}{n-sp}=p_s^*,$$ 
and again we conclude that 
$$H^{s,p}(\R^n)\dhookrightarrow{} L^q(\Omega),$$ 
for all $1\leq q<p_s^*$. Finally, we assume that $sp<n$ and $p>n$. We have that $W^{1,p}(\R^n)\dhookrightarrow L^q(\Omega)$ for $q\in [1,+\infty]$. By the same argument as above, 
$$H^{s,p}(\R^n)\dhookrightarrow L^{r_s(q)}(\Omega),$$ 
for any $r_s(q)=\left(\frac{1-s}{p}+\frac{s}{q}\right)^{-1}$, with $q\in [p,+\infty]$. Since $r_s(p)=p$ and $r_s(+\infty)=\frac{p}{1-s}$, we conclude that 
$$H^{s,p}(\R^n)\dhookrightarrow L^q(\Omega),$$ 
for any $1\leq q\leq \frac{p}{1-s}$. Now, let $q\in \left(p/(1-s),p_s^*\right)$. The complex interpolation method satisfies that for every Banach space $E$, $[E,E]_\theta=E$, $\theta\in (0,1)$, with the same norm. Hence, choosing 
$$\alpha=\frac{n\left(q(1-s)-p\right)}{sq(p-n)}\in (0,1),$$
we have that 
$$H^{s,p}(\R^n)=[H^{s,p}(\R^n),H^{s,p}(\R^n)]_\alpha\dhookrightarrow [L^{p/(1-s)}(\Omega),L^{p_s^*}(\Omega)]_\alpha=L^q(\Omega),$$ therefore $$H^{s,p}(\R^n)\dhookrightarrow L^q(\Omega),$$
for all $1\leq q<p_s^*$, which concludes the proof of $a)$. 

%\begin{Teor}[Fractional Rellich-Kondrachov-Critical case]
    %Let $\Omega\subset \R^n$ a bounded set, $s\in (0,1)$ and $1<p<\infty$ such that $sp=n$. Then, $$H^{s,p}(\R^n)\xhookrightarrow{}\xhookrightarrow{}L^q(\Omega),\,1<q<\infty.$$
%\end{Teor}
%\noindent{}\textbf{Proof:} 
\item[b)] Let $\varepsilon\in (0,s)$. By {Theorem \ref{BesselProps}}, 
$$H^{s,p}(\R^n)\xhookrightarrow{}H^{s-\varepsilon,p}(\R^n),$$ 
and since $(s-\varepsilon)p<n$, by $a)$, 
$$H^{s-\varepsilon,p}(\R^n)\dhookrightarrow L^q(\Omega),\,1<q<p_{s-\varepsilon}^*,$$ hence by \textbf{Corollary \ref{compactcompos}}, 
$$H^{s,p}(\R^n)\dhookrightarrow L^q(\Omega),$$
for all $1\leq q<p_{s-\varepsilon}^*$. Since $\varepsilon$ is arbitrary and the function $\varepsilon\mapsto p_{s-\varepsilon}^*=\frac{np}{n-(s-\varepsilon)p}$ is a continuous function mapping $(0,s)$ into $(p,\infty)$, we conclude that $$H^{s,p}(\R^n)\dhookrightarrow L^q(\Omega), $$
for all $1\leq q<\infty$.

%\begin{Teor}[Fractional Rellich-Kondrachov-Supercritical case]
 %   Let $\Omega\subset\R^n$ with finite measure. Let $s\in (0,1)$ and $1<p<\infty$ such that $sp>n$. Then, $$H^{s,p}(\R^n)\xhookrightarrow{}\xhookrightarrow{}C^{0,\mu}(\overline{\Omega}),\,0<\mu<\mu_s=s-n/p.$$
%\end{Teor}
%\noindent{}\textbf{Proof:} 

\item[c)] Since $sp>n$, $p>n$, and by the Rellich-Kondrachov theorem for classical Sobolev spaces, 
$$W^{1,p}(\R^n)\dhookrightarrow C^{0,\mu}(\overline{\Omega}),$$
for $0<\mu<\mu_1=1-n/p$. Let $0<\mu<\mu_1$ arbitrary and consider the inclusion operator $$i:L^p(\R^n)+W^{1,p}(\R^n)\to L^p(\overline{\Omega})+C^{0,\mu}(\overline{\Omega}),$$ which is and admissible operator such that $i:W^{1,p}(\R^n)\to C^{0,\mu}(\overline{\Omega})$ is compact. Then, since $W^{1,p}(\R^n)$ is an UMD-space, $$i:[L^p(\R^n),W^{1,p}(\R^n)]_s\to [L^p(\Omega),C^{0,\mu}(\overline{\Omega})]_s,$$ is a compact inclusion by {Theorem \ref{ComplexCompact}}. We have that $[L^p(\R^n),W^{1,p}(\R^n)]_s=H^{s,p}(\R^n)$ and $[L^p(\Omega),C^{0,\mu}(\overline{\Omega})]_s$ is embedded into the Morrey-Campanato space $\mathcal{L}^{p,s(p\mu+n)}(\Omega)$ (see \cite{BellidoGarciaA2025}). If we restrict ourserveles to $\mu\in \left(\frac{n(1-s)}{sp},\mu_1\right)$, we have that $$s(p\mu+n)=sp\mu+sn>sp\frac{n(1-s)}{sp}+sn=n-sn+sn=n,$$ and hence $\mathcal{L}^{p,s(p\mu+n)}(\Omega)=C^{0,\alpha_s(\mu)}(\overline{\Omega})$, where $\alpha_s(\mu)=\frac{sp\mu+sn-n}{p}$. The function $\mu\mapsto \alpha_s(\mu)$ is a continuous function such that $$\lim_{\mu\to \left(\frac{n(1-s)}{sp}\right)^+}\alpha_s(\mu)=0,$$ and $$\lim_{\mu\to \mu_1^-}\alpha_s(\mu)=\frac{sp(1-n/p)+sn-n}{p}=\frac{sp-n}{p}=s-\frac{n}{p}=\mu_s,$$ and hence it maps continuously $\left(\frac{n(1-s)}{sp},\mu_1\right)$ into $(0,\mu_s)$, and hence 
$$H^{s,p}(\R^n)\dhookrightarrow C^{0,\alpha}(\overline{\Omega}),\,0<\alpha<\mu_s.$$
\end{itemize}\qed 

Finally, the following proposition gives a compactness result between Bessel potential spaces. 

\begin{Prop}
    Let $0<t<s<1$, $1<p<\infty$ and $\Omega\subset \R^n$ be a bounded subset with Lipschitz boundary. Then, 
    $$H^{s,p}(\R^n)\dhookrightarrow H^{t,p}(\Omega).$$
\end{Prop}
\noindent{}\textbf{Proof:} By the classical Rellich-Kondrachov Theorem for Sobolev spaces, $W^{1,p}(\R^n)\dhookrightarrow L^p(\Omega)$, and hence by \textbf{Theorem \ref{functorialPeetreIII}}, the embedding $$H^{s,p}(\Omega)=[L^p(\Omega),W^{1,p}(\Omega)]_s\xhookrightarrow{}[L^p(\Omega),W^{1,p}(\Omega)]_t=H^{t,p}(\Omega),$$ is compact, so by \textbf{Corollary \ref{compactcompos}}, $$H^{s,p}(\R^n)\dhookrightarrow H^{t,p}(\Omega).\qed$$

\subsection{A proof based on estimations on translations}

%%%%%%%%%%%%%%%%%%%%%%%%%%%%%%%%%%%%%%%%%%%%%%%%%%%%%%%%%%%%%%%%%%%%%

In this section we prove Theorem \ref{BesselCompact} as a consequence of Fréchet-Kolmogorov-Riesz theorem, and for it we need to show estimates on the translation-difference operator $\tau_h^-$, with $$\tau_h^-:u\to u(\cdot+h)-u,\,u\in C^\infty_c(\R^n),\,h\in \R^n.$$ To this aim, we consider the identification of the Bessel potential space $H^{s,p}(\R^n)$ with $\overline{C_c^\infty(\R^n)}^{\norm{\cdot}_{D^s,p}}$. In particular, we are going to establish that for every $s\in (0,1)$ and $1\leq p<\infty$, there exists a positive constant $C$, not depending on $s$, such that 
$$\norm{\tau_h^-u}_p\leq \frac{C}{s(1-s)}|h|^s\norm{D^su}_p,\,u\in H^{s,p}(\R^n).$$ 
This inequality generalizes the classical result for the Sobolev space $W^{1,p}(\R^n)$ \cite[Proposition~9.3 (iii)]{Brezis2010}.
The case $p=1$ was proved in \cite[Proposition~3.14]{ComiStefani2019}, and here we generalize the estimation for $p\ge 1$. It is also worth mentioning that a similar idea has been used recently for proving compact embeddings for Gagliardo spaces in \cite{TesoGomezVazquez2020}, proving estimations for the translation-difference operator based on the Gagliardo seminorm using fine properties of the $K$-functional of the real interpolation method. Interestingly, this approach allows us to cover case $p=1$ in Theorem \ref{BesselCompact}, which cannot be covered with interpolation as equivalence of $H^{s,p}(\R^n)$ with $[L^p(\R^n),W^{1,p}(\R^n)]_s$ only holds for $1<p<\infty$ \cite{BellidoGarciaA2025}. See also this reference for a discussion on the problem of whether or not $H^{s,1}(\R^n)$ is an interpolation space. The following is a technical lemma.
\begin{Lema}\label{CuetoLema}
    Let $s\in (0,1)$ and the function $u(x)=x/|x|^{n+1-s}$. Then, there exists a positive constant $C$, independent of $s$, such that $$\norm{\tau_{-e_1}^-u}_1\leq \frac{C}{s(1-s)},$$ where $e_1$ is the first vector of the canonical basis of $\R^n$.
\end{Lema}
\noindent{}\textbf{Proof:} First, we see that \begin{align*}
    \int_{B(0,2)}|\tau_{-e_1}^-u(x)\,dx|&=\int_{B(0,2)}\left|\frac{x}{|x|^{n+1-s}}-\frac{x-e_1}{|x-e_1|^{n+1-s}}\right|\,dx\leq\int_{B(0,2)}\left(\frac{1}{|x|^{n-s}}+\frac{1}{|x-e_1|^{n-s}}\right)\,dx\\
    &\leq C\int_{B(0,2)}\frac{1}{|x|^{n-s}}\,dx=C\omega_n\frac{2^s}{s}\leq \frac{2C\omega_n}{s},
\end{align*} where $C>0$ does not depend on $s$. Now, for a fixed $x\in B(0,2)^c$, \begin{align*}
    |\tau_{-e_1}^-u(x)|&=\left|\int_0^1\frac{d}{dt}\left(\frac{x-te_1}{|x-te_1|^{n+1-s}}\right)\,dt\right|\\
    &=\left|\int_0^1 (n+1-s)\frac{[(x-te_1)\cdot e_1](x-te_1)}{|x-te_1|^{n+3-s}}-\frac{e_1}{|x-te_1|^{n+1-s}}\,dt\right|\\
    &\leq \int_0^1\left(\frac{n+1}{|x-te_1|^{n+1-s}}+\frac{1}{|x-te_1|^{n+1-s}}\right)\,dt=(n+2)\int_0^1\frac{1}{|x-te_1|^{n+1-s}}\,dt.
\end{align*}
Since $|x|>2$ and $0\leq t\leq 1$, $$|x-te_1|\geq \left||x|-|te_1|\right|=\left|x|-t\right|=|x|-t\geq |x|-1\geq \frac{|x|}{2},$$ thus $$\int_0^1\frac{1}{|x-te_1|^{n+1-s}}\,dt\leq \int_0^1\frac{2^{n+1-s}}{|x|^{n+1-s}}\,dt=\frac{2^{n+1-s}}{|x|^{n+1-s}}\leq \frac{2^{n+1}}{|x|^{n+1-s}}.$$ Then, \begin{align*}
    \int_{B(0,2)^c}|\tau_{-e_1}^-u(x)|\,dx&\leq 2^{n+1}(n+2)\int_{B(0,2)^c}\frac{1}{|x|^{n+1-s}}\,=2^{n+1}(n+2)\omega_n\int_2^\infty r^{-2+s}\,dr\\
    &=2^{n+1}(n+2)\omega_n\frac{2^{s-1}}{1-s}=\frac{\omega_n(n+2)2^{n+s}}{1-s}\leq \frac{\omega_n(n+2)2^{n+1}}{1-s}.
\end{align*}
Putting all together we get 
\begin{align*}\norm{\tau_{-e_1}^-u(x)}_1&=\int_{B(0,2)}|\tau_{-e_1}^-u(x)|\,dx+\int_{B(0,2)^c}|\tau_{-e_1}^-u(x)|\,dx\leq \frac{2C\omega_n}{s}+\frac{\omega_n(n+2)2^{n+s}}{1-s}\\
&\le \frac{C'}{s(1-s)},\end{align*} where $C'>0$ does not depend on $s$.\qed\\

We are finally ready to prove estimation on the the translation-difference operator based on the Riesz fractional gradient. 
\begin{Prop}[Estimation on translations for Bessel potential spaces]\label{BesselEstimate}
    Let $s\in (0,1)$ and $1\leq p<\infty$. Then, there exists a positive constant $C$, independent of $s$ and $p$, such that $$\norm{\tau_h^-u}_p\leq \frac{C}{s(1-s)}|h|^s\norm{D^su}_p,\,u\in H^{s,p}(\R^n),$$ for every $h\in \R^n$.
\end{Prop}
\noindent{}\textbf{Proof:} Since $H^{s,p}(\R^n)=\overline{C_c^\infty(\R^n)}^{\norm{\cdot}_{D^s,p}}$, it is enough to prove it for $u\in C_c^\infty(\R^n)$. Let $h\in \R^n$ fixed. By Proposition \ref{FTOC}, $$u(x+h)-u(x)=c_{n,-s}\int_{\R^n}\left(\frac{x+h-y}{|x+h-y|^{n+1-s}}-\frac{x-y}{|x-y|^{n+1-s}}\right)\cdot D^su(y)\,dy,$$ 
and making the change of variables $y-x=z$, \begin{align*}
    u(x+h)-u(x)&=c_{n,-s}\int_{\R^n}\left(\frac{h-z}{|h-z|^{n+1-s}}-\frac{-z}{|-z|^{n+1-s}}\right)\cdot D^su(x+z)\,dz\\
    &=c_{n,-s}\int_{\R^n}\left(\frac{z}{|z|^{n+1-s}}-\frac{z-h}{|z-h|^{n+1-s}}\right)\cdot D^su(x+z)\,dz,
\end{align*}
 so \begin{align*}
     |u(x+h)-u(x)&=|c_{n,-s}|\left|\int_{\R^n}\left(\frac{z}{|z|^{n+1-s}}-\frac{z-h}{|z-h|^{n+1-s}}\right)\cdot D^su(x+z)\,dz\right|\\
     &\leq |c_{n,-s}|\int_{\R^n}\left|\frac{z}{|z|^{n+1-s}}-\frac{z-h}{|z-h|^{n+1-s}}\right||D^su(x+z)|\,dz.
 \end{align*}
Let $R\in O(n)$ (the set of proper rotations of $\R^n$), with $\det(R)=-1$ such that $R^Th=|h|e_1$, the first vector of the canonical basis of $\R^n$. Then, making the change of variables $z=|h|Rw$, since $R$ has determinant equal to $-1$ and $z-h=|h|Rw-R|h|e_1=|h|R(w-e_1)$,  \begin{align*}&\int_{\R^n}\left|\frac{z}{|z|^{n+1-s}}-\frac{z-h}{|z-h|^{n+1-s}}\right||D^su(x+z)|\,dz\\
&=|h|^s\int_{\R^n}\left|\frac{w}{|w|^{n+1-s}}-\frac{w-e_1}{|w-e_1|^{n+1-s}}\right||D^s\left(x+|h|Rw\right)|\,dw\\
&=|h|^s\int_{\R^n}\left|\frac{w}{|w|^{n+1-s}}-\frac{w-e_1}{|w-e_1|^{n+1-s}}\right|^{1/p+1/q}|D^s\left(x+|h|Rw\right)|\,dw\end{align*} were $1p+1/q=1$. Now, by H\"older's inequality 
\begin{align*}&\int_{\R^n}\left|\frac{w}{|w|^{n+1-s}}-\frac{w-e_1}{|w-e_1|^{n+1-s}}\right|^{1/p+1/q}|D^su(x+|h|Rw)|\,dw\\
&\leq \left(\int_{\R^n}\left|\frac{w}{|w|^{n+1-s}}-\frac{w-e_1}{|w-e_1|^{n+1-s}}\right|^{q/q}\,dw\right)^{1/q}\\
&\left(\int_{\R^n}\left|\frac{w}{|w|^{n+1-s}}-\frac{w-e_1}{|w-e_1|^{n+1-s}}\right|^{p/p}|D^su(x+|h|Rw)|^p\,dw\right)^{1/p}
\end{align*} and by {Lemma \ref{CuetoLema}}, \begin{align*}&\int_{\R^n}\left|\frac{w}{|w|^{n+1-s}}-\frac{w-e_1}{|w-e_1|^{n+1-s}}\right|^{1/p+1/q}|D^su(x+|h|Rw)|\,dw\\
&\leq \left(\frac{C}{s(1-s)}\right)^{1/q}\left(\int_{\R^n}\left|\frac{w}{|w|^{n+1-s}}-\frac{w-e_1}{|w-e_1|^{n+1-s}}\right||D^su(x+|h|Rw)|^p\,dw\right)^{1/p}.
\end{align*}
Putting all together and raising to the $p$ yields \begin{align*}
    &|u(x+h)-u(x)|^p \\
    &\leq |c_{n,-s}|^p|h|^{sp}\left(\frac{C}{s(1-s)}\right)^{p/q}\int_{\R^n}\left|\frac{w}{|w|^{n+1-s}}-\frac{w-e_1}{|w-e_1|^{n+1-s}}\right||D^su(x+|h|Rw)|^p\,dw
\end{align*}
so using Fubini's theorem and {Lemma \ref{CuetoLema}} again, \begin{align*}
    &\norm{\tau_h^-u}_p^p=\int_{\R^n}|u(x+h)-u(x)|^p\,dx\\
    &\leq |c_{n,-s}|^p|h|^{sp}\left(\frac{C}{s(1-s)}\right)^{p/q}\int_{\R^n}\left(\int_{\R^n}\left|\frac{w}{|w|^{n+1-s}}-\frac{w-e_1}{|w-e_1|^{n+1-s}}\right||D^su(x+|h|Rw)|^p\,dw\right)\,dx\\
    &=|c_{n,-s}|^p|h|^{sp}\left(\frac{C}{s(1-s)}\right)^{p/q}\int_{\R^n}\left|\frac{w}{|w|^{n+1-s}}-\frac{w-e_1}{|w-e_1|^{n+1-s}}\right|\left(\int_{\R^n}|D^su(x+|h|Rw)|^p\,dx\right)\,dw\\
    &= |c_{n,-s}|^p|h|^{sp}\left(\frac{C}{s(1-s)}\right)^{p/q}\norm{D^su}_p^p \int_{\R^n}\left|\frac{w}{|w|^{n+1-s}}-\frac{w-e_1}{|w-e_1|^{n+1-s}}\right|\,dw\\
    &\leq |c_{n,-s}|^p|h|^{sp}\left(\frac{C}{s(1-s)}\right)^{p/q+1}\norm{D^su}_p^p= |c_{n,-s}|^p|h|^{sp}\left(\frac{C}{s(1-s)}\right)^{p}\norm{D^su}_p^p,
\end{align*}
so $$\norm{\tau_h^-u}_p\leq |c_{n-s}||h|^s\frac{C}{s(1-s)}\norm{D^su}_p\leq \frac{\Tilde{C}}{s(1-s)}|h|^s\norm{D^su}_p,$$ where $$\Tilde{C}=C\operatorname{sup}_{s\in (0,1)}c_{n,-s}<\infty,$$ and does not depend on $s$.\qed\\

\textbf{Proof of Theorem \ref{BesselCompact}:}\begin{itemize}
 \item[a)] Let $H\subset H^{s,p}(\R^n)$ an uniformly bounded subset, and hence uniformly bounded on $L^p(\Omega)$, and denote $M:=\operatorname{sup}_{u\in H}\{\norm{u}_{D^s,p}\}$.
By {Theorem \ref{BesselEstimate}}, $$\norm{\tau_h^-u}_p\leq \frac{C}{s(1-s)}|h|^s\norm{D^su}_{p}\leq \frac{C}{s(1-s)}|h|^s\norm{u}_{D^{s},p}\leq \frac{MC}{s(1-s)}|h|^s\to 0,\,\text{as}\,|h|\to 0^+,\,u\in H.$$ Given $\varepsilon>0$ arbitrary, we can choose $\delta=\left(\frac{s(1-s)\varepsilon}{MC}\right)^{1/s}$, which not depends on the function $u$, so $\norm{\tau_h^-u}_p<\varepsilon$ whenever $|h|<\delta$ uniformly on $u\in W$. Hence, by {Theorem \ref{FrechetKolmogorov}}, $H^{s,p}(\R^n)\dhookrightarrow L^p(\Omega)$, and since $\Omega$ is bounded, $H^{s,p}(\R^n)\dhookrightarrow L^q(\Omega)$ for every $1\leq q\leq p$. Now, let $q\in (p,p_s^*)$ and $\{u_m\}$ a bounded sequence in $H^{s,p}(\R^n)$, and hence also bounded in $L^p(\R^n)$. Since $H^{s,p}(\R^n)\dhookrightarrow L^p(\Omega)$, by the sequence characterization of compact embddings, there exists a subsequence $\{u_{m_k}\}$ such that is Cauchy for the $L^p(\Omega)$-norm. Since 
$$[L^p(\Omega),L^{p_s^*}(\Omega)]_\alpha=L^q(\Omega),\,\alpha=\frac{n(q-p)}{sqp}\in (0,1),$$ 
we have that for any $k,l\in\mathbb{N}$, 
$$\norm{u_{m_k}-u_{m_l}}_{L^q(\Omega)}\leq \norm{u_{m_k}-u_{m_l}}_{L^p(\Omega)}^{1-\alpha}\norm{u_{m_k}-u_{m_l}}_{L^{p_s^*}(\Omega)}^{\alpha}\to 0,$$ 
as $k,l\to \infty$, due to $\{u_{m_k}\}$ being Cauchy in $L^p(\Omega)$ and the continuous embedding $H^{s,p}(\R^n)\xhookrightarrow{}L^{p_s^*}(\R^n)$. We conclude then that $\{u_{m_k}\}$ is a Cauchy subsequence of $\{u_m\}$ in $L^q(\Omega)$, so $$H^{s,p}(\R^n)\dhookrightarrow L^q(\Omega),\,1\leq q<p_s^*.$$

\item[b)] With the same argument as in the previous item we show that $$H^{s,p}(\R^n)\dhookrightarrow{} L^q(\Omega),\,1\leq q\leq p,$$ and that for every bounded sequence $\{v_m\}\subset H^{s,p}(\R^n)$, there exists a Cauchy subsequence $\{v_{m_k}\}$ in the $L^p(\Omega)$-norm. Now, take $q\in (p,\infty)$ and $q<r<\infty$. Since $$[L^p(\Omega),L^r(\Omega)]_\beta=L^q(\Omega),\,\beta=\frac{r(q-p)}{q(r-p)}\in (0,1),$$ we have that for any $k,l\in\mathbb{N}$, $$\norm{v_{m_k}-v_{m_l}}_{L^q(\Omega)}\leq \norm{v_{m_k}-v_{m_l}}_{L^p(\Omega)}^{1-\beta}\norm{v_{m_k}-v_{m_l}}_{L^{r}(\Omega)}^{\beta}\to 0,$$ as $k,l\to \infty$, due to $\{u_{m_k}\}$ being Cauchy in $L^p(\Omega)$ and the continuous embedding $H^{s,p}(\R^n)\xhookrightarrow{}L^{r}(\R^n)$. We conclude then that $\{v_{m_k}\}$ is a Cauchy subsequence of $\{v_m\}$ in $L^q(\Omega)$, so $$H^{s,p}(\R^n)\dhookrightarrow L^q(\Omega),\,1\leq q<\infty.$$

\item[c)] For $sp>n$ we have that \cite[Theorem 3.21 (3)]{BellidoGarciaA2025} $$H^{s,p}(\R^n)\xhookrightarrow{}C^{0,\mu_s^*}(\R^n).$$ 
Let $0<\beta<1$ and consider $\{f_m\}\subset C^{0,\beta}(\overline{\Omega})$ an uniformly bounded sequence. Set \\
$F:=\operatorname{sup}_{m\in\mathbb{N}}\{\norm{f_m}_{C^{0,\beta}(\overline{\Omega})}\}$. Since for every $m\in \mathbb{N}$ and $x,y\in \overline{\Omega}$ we have that $$|f_n(x)-f_n(y)|\leq F|x-y|^\beta,$$ for any $\varepsilon>0$ we can choose $\delta=(\delta/F)^{1/\beta}$ and then for every $x,y\in \overline{\Omega}$ such that $|x-y|<\delta$, we have $$|f_m(x)-f_m(y)|<\varepsilon,\,\forall m\in\mathbb{N},$$ thus by {Theorem \ref{ascoliarzela}}, there exists an uniformly convergent subsequence $\{f_{m_k}\}\subset \{f_m\}$. Now, let $0<\alpha<\beta$. For every $k\in\mathbb{N}$ and $x,y\in \overline{\Omega}$ we have that $$\frac{|f_{m_k}(x)-f_{m_k}(y)|}{|x-y|^\alpha}=\left(\frac{|f_{m_k}(x)-f_{m_k}(y)|}{|x-y|^\beta}\right)^{\alpha/\beta}|f_{m_k}(x)-f_{m_k}(y)|^{1-\alpha/\beta}\leq 2F[f_{m_k}]^{\alpha/\beta}_{C^{0,\beta}(\overline{\Omega})},$$ and hence $$\norm{f_{m_k}}_{C^{0,\alpha}(\overline{\Omega})}\leq 2F\norm{f_{m_k}}_{C^{0,\beta}(\overline{\Omega})},$$ so $\{f_{m_k}\}$ is a Cauchy subsequence in $C^{0,\alpha}(\overline{\Omega)}$, which in particular implies that 
for every $0<\beta<\alpha<1$, $$C^{0,\beta}(\overline{\Omega})\dhookrightarrow C^{0,\alpha}(\overline{\Omega}).$$ Then, by {Corollary \ref{compactcompos}}, $$H^{s,p}(\R^n)\dhookrightarrow{} C^{0,\mu}(\overline{\Omega}),\,0<\mu<\mu_s^*.$$

%\item[d)] As in the cases a) and b) we have the compact embedding $H^{s,p}(\R^n)\xhookrightarrow{}\xhookrightarrow{}L^p(\Omega)$, and hence for every bounded sequence $\{w_m\}\subset H^{s,p}(\R^n)$, there exists a subsequence $\{w_{m_k}\}$ that is Cauchy in $L^p(\Omega)$. Since its bounded in $H^{s,p}(\R^n)$, it is also bounded in $H^{s,p}(\Omega)$, and since $\Omega$ has Lipschitz domain $$[L^p(\Omega),H^{s,p}(\Omega)]_{t/s}=H^{t,p}(\Omega),$$ so there exists a positive constant depending continuously on the parameters $n,p,s$ such that $$\norm{w_{m_k}-w_{m_l}}_{H^{t,p}(\Omega)}\leq C\norm{w_{m_k}-w_{m_l}}_{L^p(\Omega)}^{1-t/s}\norm{w_{m_k}-w_{m_l}}_{H^{s,p}(\Omega)}^{t/s}\to 0,$$ as $k,l\to \infty$. Hence, there exists a subsequence $\{w_{m_k}\}\subset \{w_m\}$ such that is Cauchy for $H^{t,p}(\Omega)$, so $$H^{s,p}(\R^n)\xhookrightarrow{}\xhookrightarrow{}H^{t,p}(\Omega).$$

\end{itemize}\qed

\subsection{Compactness and contiguity}

Nested embeddings between Gagliardo and Bessel spaces provided by Theorem \ref{Contiguity} make both families of spaces intimately related. Actually, from Theorem \ref{Contiguity}, it is easy to prove that embeddings, both continuous and compact, for Gagliardo and Bessel potential spaces are the same. As a example, we prove here Theorem \ref{BesselCompact}, in the case $1<p<\infty$, as a consequence of the following compact embeddings for Gagliardo spaces. Proof of the following result can be found in \cite{DiNezzaPalatucciValdinoci2012,TesoGomezVazquez2020}. 

\begin{Teor}[Rellich-Kondrachov for Gagliardo spaces]\label{GagliardoCompact}
    Let $1\le p<\infty$, $s\in (0,1)$ and $\Omega\subset \R^n$ a measurable subset with finite measure. Then, the following holds:
    \begin{itemize}
        \item[a)] If $sp<n$, then $W^{s,p}(\R^n)\dhookrightarrow L^q(\Omega)$ for $1\leq q<p_s^*=\frac{np}{n-sp}$.
        \item[b)] If $sp=n$, then $W^{s,p}(\R^n)\dhookrightarrow L^q(\Omega)$ for $1\leq q<\infty$.
        \item[c)] If $sp>n$, then $W^{s,p}(\R^n)\dhookrightarrow C^{0,\mu}(\overline{\Omega})$ for $0<\mu<\mu_s^*=s-n/p$.
        %\item[d)]If $s>t$ and $\Omega$ is bounded with Lipschitz boundaru, then $W^{s,p}(\R^n)\xhookrightarrow{}\xhookrightarrow{}W^{t,p}(\Omega)$.
    \end{itemize}
\end{Teor}
%The proof will be given in the following subsection. Now we prove the theorem for Bessel spaces.
%\begin{Teor}[Rellich-Kondrachov for Bessel potential spaces]
 %   Let $s\in (0,1)$, $1<p<\infty$ and $\Omega\subset \R^n$ a bounded subset. Then, \begin{itemize}
  %      \item[a)] If $sp<n$, then $H^{s,p}(\R^n)\xhookrightarrow{}\xhookrightarrow{}L^q(\Omega)$ for $1\leq q<p_s^*=\frac{np}{n-sp}$.
   %     \item[b)] If $sp=n$, then $H^{s,p}(\R^n)\xhookrightarrow{}\xhookrightarrow{}L^q(\Omega)$ for $1\leq q<\infty$.
   %     \item[c)] If $sp>n$, then $H^{s,p}(\R^n)\xhookrightarrow{}\xhookrightarrow{}C^{0,\mu}(\overline{\Omega})$ for $0<\mu<\mu_s^*=s-n/p$.
   % \end{itemize}
%\end{Teor}
\noindent\textbf{Proof of Theorem \ref{BesselCompact} (case $1<p<\infty$):} 

\begin{itemize}
    \item[a)] If $sp<n$, by {Theorem \ref{GagliardoCompact} a)}, $W^{s,p}(\R^n)\dhookrightarrow L^q(\Omega)$ for $1<q<p_s^*$. Now, let $\varepsilon\in (0,s)$. By  {Theorem \ref{Contiguity}}, 
    $$H^{s,p}(\R^n)\xhookrightarrow{}W^{s-\varepsilon,p}(\R^n)\dhookrightarrow L^q(\Omega),$$ for $1<q<p_{s-\varepsilon}^*=\frac{np}{n-(s-\varepsilon)p}$, and hence by {Proposition \ref{compactcompos}}, 
    $$H^{s,p}(\R^n)\dhookrightarrow L^q(\Omega),\,1<q<p_{s-\varepsilon}^*.$$ 
    The function $\varepsilon\mapsto p_{s-\varepsilon}^*$ is a continuous decreasing function such that 
    $$\lim_{\varepsilon\to 0^+}p_{s-\varepsilon}^*=\lim_{\varepsilon\to 0^+}\frac{np}{n-(s-\varepsilon)p}=p_s^*,$$ and 
    $$\lim_{\varepsilon\to s^-}p_{s-\varepsilon}^*=\lim_{\varepsilon\to s^-}\frac{np}{n-(s-\varepsilon)p}=p,$$ 
    so it maps continuously $(0,s)$ into $(p,p_s^*)$, and since $\varepsilon$ is arbitrary we get that $$H^{s,p}(\R^n)\dhookrightarrow L^q(\Omega),$$ for $q\in (p,p_s^*)$. Now, since $\Omega$ is bounded, we conclude that 
    $$H^{s,p}(\R^n)\dhookrightarrow L^q(\Omega),$$
    for all $q\in (1,p_s^*)$.
    
    \item[b)] If $sp=n$, with an analogous argument we get that 
    $$H^{s,p}(\R^n)\dhookrightarrow L^q(\Omega),\,q\in (1,p_{s-\varepsilon}^*).$$ Again, the function $\varepsilon\mapsto p_{s-\varepsilon}^*$ is a continuous decreasing function such that 
    $$\lim_{\varepsilon\to 0^+}p_{s-\varepsilon}^*=\lim_{\varepsilon\to 0^+}\frac{np}{n-(s-\varepsilon)p}=+\infty,$$ 
    and 
    $$\lim_{\varepsilon\to s^-}p_{s-\varepsilon}^*=\lim_{\varepsilon\to s^-}\frac{np}{n-(s-\varepsilon)p}=p,$$ 
    so it maps continuously $(0,s)$ into $(p,+\infty)$, and since $\varepsilon$ is arbitrary we get that 
    $$H^{s,p}(\R^n)\dhookrightarrow L^q(\Omega),$$ 
    for 
    $q\in (p,+\infty)$. Now, since $\Omega$ is bounded, we conclude that 
    $$H^{s,p}(\R^n)\dhookrightarrow L^q(\Omega),$$
    for all $q\in (1,+\infty)$.
    
    \item[c)] If $sp>n$, by {Theorem \ref{GagliardoCompact} c)}, we have that $W^{s,p}(\R^n)\dhookrightarrow C^{0,\mu}(\overline{\Omega})$ for $0<\mu<\mu_s^*$. Now, by {Theorem \ref{Contiguity}} we have that 
    $$H^{s,p}(\R^n)\xhookrightarrow{}W^{s-\varepsilon,p}(\R^n),$$ 
    for $\varepsilon\in (0,\mu_s)$, and hence 
    $$H^{s,p}(\R^n)\xhookrightarrow{}W^{s-\varepsilon,p}(\R^n)\dhookrightarrow C^{0,\mu}(\overline{\Omega}),$$ 
    for $0<\mu<\mu_{s-\varepsilon}^*=(s-\varepsilon)-n/p$. The function $\varepsilon\mapsto \mu_{s-\varepsilon}^*$ maps continuously $(0,\mu_s)$ into $(0,\mu_s^*)$, and hence by {Proposition \ref{compactcompos}}, 
    $$H^{s,p}(\R^n)\dhookrightarrow C^{0,\mu}(\overline{\Omega}),\,0<\mu<\mu_s^*.$$ \qed
\end{itemize}

\section*{Acknowledgements}

The work of JCB and GG-Shas been supported by {\it Agencia Estatal de Investigación} (Spain) through grant PID2023-151823NB-I00 and {\it Junta de Comunidades de Castilla-La Mancha} (Spain) through grant SBPLY/23/180225/000023. GG-S is funded by a Doctoral Fellowship by Universidad de Castilla-La Mancha 2024-UNIVERS-12844-404. The work of JC has been supported by the {\it Agencia Estatal de Investigación} (Spain) through the grant PID2021-124195NB-C32, and by the Madrid Government (Comunidad de Madrid, Spain) under the multiannual Agreement with UAM in the line for the Excellence of the University Research Staff in the context of the V PRICIT (Regional Programme of Research and Technological Innovation).

\section*{Conflicts of interest}

The authors declare that there are no conflicts of interest regarding the publication of this paper.
%-------------------------------------------
%-------------------------------------------
%-------------------------------------------
% Bibliography
%-------------------------------------------
%-------------------------------------------
%-------------------------------------------
%\bibliography{References}{}
%\bibliographystyle{plain} 

\end{document}